\documentclass[10pt,draftclsnofoot,onecolumn]{IEEEtran}
\ifCLASSINFOpdf
\else
\fi
\usepackage{graphicx}
\usepackage{caption}
\usepackage{subcaption}
\usepackage{amsmath} 
\usepackage{amssymb}  
\usepackage{amsthm}
\usepackage{bm}
\usepackage{lipsum}
\usepackage[linesnumbered, ruled]{algorithm2e}
\usepackage{enumerate}

\usepackage{color}

\newtheorem{theorem}{Theorem}
\newtheorem{proposition}{Proposition}
\newtheorem{lemma}{Lemma}
\newtheorem{cor}{Corollary}

\newtheorem{example}{Example}

\begin{document}
%
\title{Risk-Limiting Dynamic Contracts for \\Direct Load Control}
%
%
%

\author{Insoon~Yang,~
Duncan S.~Callaway,~ 
and
        Claire J.~Tomlin
        \thanks{This work was supported by the NSF CPS project ActionWebs under grant number 0931843, NSF CPS project FORCES under grant number 1239166,
Robert Bosch LLC through its Bosch Energy Research Network funding program and NSF CPS Award number 1239467.}
\thanks{I. Yang and C. J. Tomlin are with the Department of Electrical Engineering and Computer Sciences, University of California, Berkeley, CA 94720, USA
        {\tt\small \{iyang, tomlin\}@eecs.berkeley.edu}}%
\thanks{D. S. Callaway is with the Energy and Resources Group, University of California, Berkeley, CA 94720, USA
        {\tt\small dcal@berkeley.edu}}
}

\maketitle

\begin{abstract}
This paper proposes a novel continuous-time dynamic contract framework that has a risk-limiting capability.
If a principal and an agent enter into such a contract, the principal can optimally manage its performance and risk with a guarantee that the agent's risk is less than or equal to a pre-specified level and that the agent's expected payoff is greater than or equal to another pre-specified threshold.
We achieve such risk-management capabilities by formulating the contract design problem as mean-variance constrained risk-sensitive control.
A dynamic programming-based method is developed to solve the problem. The key idea of our proposed solution method is to reformulate the inequality constraints on the mean and the variance of the agent's payoff as dynamical system constraints by introducing new state and control variables. The reformulations use the martingale representation theorem.
The proposed contract method enables us to develop a new direct load control method that provides the load-serving entity with financial risk management solutions in real-time electricity markets.
We also propose 
an approximate decomposition of the optimal contract design problem for multiple customers into multiple low-dimensional contract problems for one customer.
This allows the direct load control program to work with a large number of customers without any scalability issues. Furthermore, the contract design procedure can be completely parallelized. 
The performance and usefulness of the proposed contract method and its application to direct load control are demonstrated using data on the electric energy consumption of customers in Austin, Texas as well as the Electricity Reliability Council of Texas' locational marginal price data.
\end{abstract}

%

%
\IEEEpeerreviewmaketitle

\section{Introduction}

To reduce the greenhouse gases caused by electricity generation, 
there has been growing interest in and efforts toward integrating renewable energy sources into the electric power grid.
California, for example, plans to use renewable resources to serve 33\% of the electricity load by 2020 \cite{CPUC}.
In particular, the penetration of solar and wind energy resources is expected to significantly increase.
However, the utilization of these resources is challenging because they are uncertain and intermittent.
To absorb the uncertainty in solar and wind power, for example, the reserve capacity must be sufficiently large. 
In California, the increase in the reserve costs needs to be compensated by all customers and the amount of the compensation per customer may not be negligible when California achieves its goal of 33\% renewable energy penetration \cite{Varaiya2011}. 

Supply-side approaches have been proposed to address the uncertainty of renewable resources in economic dispatch and unit commitment using stochastic dynamic programming \cite{Rajagopal2013}, mixed-integer stochastic programming \cite{Bouffard2008} and robust optimization \cite{Bertsimas2013}, among others.
However, these aforementioned methods do not examine the potential of demand-side resources in managing uncertain renewable generation or loads. 
To investigate this potential, this paper proposes  a demand-side solution that manages the uncertainty of customers' solar and wind generation and loads. 
The load-serving entity or the aggregator for the customers procures power or generation reserves in a day-ahead market and the amount of procurement is determined based on a day-ahead load forecast.
Because the actual total load deviates from the procured power, the load-serving entity must purchase the deviated amount of power in a real-time market to balance supply and demand.
It is desirable for the load-serving entity to minimize the real-time purchase of energy because the energy price and the reserve cost in the real-time market are normally higher and more volatile than those in the day-ahead market.
As the penetration of customers' solar and wind generation increases, however,
the electricity price in the real-time market is highly volatile and the customers' demand is very difficult to predict.
Therefore,
the risk of spending a substantial budget in the real-time market increases.
If the load-serving entity bears this financial risk, the energy price in the customers' tariff would inevitably increase.
To reduce this risk that the load-serving entity must face, we propose a contract approach for \emph{direct load control} in which the customer transfers the authority to control his or her load to the load-serving entity.
Once the load-serving entity and its customer enter into the contract, the load-serving entity can allocate a portion of the risk to the customer through the compensation scheme and the control strategy for the customer's load specified in the contract.
A customer might reasonably worry that such compensation and control could increase the risk for large energy costs and a disruption of comfort.
The proposed contract addresses this concern by guaranteeing that the risk in the customer's payoff, a weighted sum of energy costs and discomfort level, is limited by a pre-specified threshold
and that the mean of the customer's payoff is greater than another pre-specified level.
The former is called the \emph{risk-limiting} condition, and the latter is called the \emph{participation payoff} condition.
The compensation scheme and the control strategy written in the contract must be designed such that their combination mitigates the load-serving entity's financial risk in the real-time market
while satisfying the risk-limiting and participation payoff conditions for the customer.

The key element of the proposed contract for such a demand-side management
 is direct load control that allows the load-serving entity to actively use the customer's load to manage its financial risk.
A number of direct control methods have been suggested for various types of electric loads such as thermostatically controlled loads \cite{Mathieu2013, Hao2013}, electric vehicles \cite{Rotering2011, Ma2013} and deferrable loads \cite{Obrien2013, Roozbehani2014}. 
The objectives of existing direct load control methods include 
shifting demand (e.g., `valley-filling'),
providing ancillary services (e.g., frequency regulation)  and energy arbitrage \cite{Callaway2011, Zhao2014, Mathieu2013a}. 
To the best of the authors' knowledge, however, the potential of direct load control for financial risk management in real-time electricity markets has not yet been studied.
We bridge the gap between direct load control and risk management by proposing a contract-based approach.

A variety of contract methods have been suggested for 
demand-side management in electric power systems. The use of contracts for reducing energy price risk in spot markets has been investigated \cite{Kaye1990}. 
Interruptible service contracts have been extensively studied, in which 
a customer takes risk of service interruption in return for a discount in the energy price \cite{Chao1987, Tan1993, Kamat2002}.
More recently, deadline-differentiated deferrable energy contracts have been proposed to prevent the risk of a customer not receiving energy delivery by a pre-specified deadline \cite{Bitar2012}.
A different contract approach associated with durations-differentiated loads is studied in \cite{Nayyar2014}.
None of the aforementioned contract methods, however, takes into account detailed electric load dynamics, which are essential in direct load control to guarantee the customer's comfort and the load's system constraints.

To incorporate dynamics of electric loads in contracts, we adopt  a dynamic contract framework, also called a continuous-time principal-agent problem \cite{Holmstrom1987, Cvitanic2012}. In such a problem, a principal (e.g., a company) and an agent (e.g., a worker) make a contract that specifies a compensation scheme and a control strategy in an uncertain environment. 
The setting we consider in this work is called the \emph{first-best}, in which the principal and the agent have the same information and share the risk in the principal's revenue stream.
The principal can monitor the agent's control or effort and, therefore, can enforce the control strategy written in the contract.
In the electricity setting, regarding a load-serving entity as the principal and its customer as the agent,
this first-best case is appropriate for direct load control because the load-serving entity has the authority to monitor and control its customer's electric loads.
A dynamic contract problem for the first-best case was first considered in \cite{Muller1998}. It uses a simple principal-agent model with exponential utility functions introduced by Holmstrom and Milgrom for moral hazard \cite{Holmstrom1987}. A more general class of the first-best case dynamic contract problems is addressed in \cite{Cadenillas2007} by using the martingale and convex duality methods \cite{Pliska1986, Cox1989, Karatzas1987, Karatzas1991}. 
However, the proposed solution approach requires that the payoff functions  be differentiable, strictly increasing and strictly concave and that the dynamical system be a stochastic integral equation. These restrictions are acceptable in many applications in economics and finance, 
but they may exclude some important engineering problems, including direct load control, because they often require dynamical systems and payoff functions be complicated.
Furthermore, the aforementioned methods assume that the utility function of the agent with respect to his or her payoff is given. 
However, in practice, it is difficult to have complete knowledge of the agent's utility function. 
In particular, if the agent's utility function used in designing the contract deviates from his or her actual utility, the agent may not want to enter into the contract again.

In this paper, we propose a novel dynamic contract method that overcomes the limitations of existing methods for the first-best case.
The proposed method uses the variance of the agent's payoff as the risk measure for the agent. By imposing a constraint on the variance, we can limit the risk the agent needs to bear. We call this constraint the agent's \emph{risk-limiting condition}.
 In addition, the proposed method guarantees that the mean of the agent's payoff exceeds some pre-specified threshold. 
From the principal's point of view, by executing an appropriately designed compensation scheme and control strategy specified in the contract,
the principal can transfer some portion of its financial risk to the agent, respecting the agent's risk-limiting condition.
This variance approach does not require complete knowledge of the agent's utility function, which is difficult to obtain in practice. 
Such a risk-limiting capability
distinguishes our method from existing contract methods.
To take into account the principal's risk aversion, we formulate the contract design problem as risk-sensitive control \cite{Jacobson1973, Fleming1995}. Due to the constraints on the mean and the variance of the agent's payoff, however, dynamic programming is not directly applicable. 
One may be able to handle the constraints using the stochastic maximum principle or the duality method \cite{Peng1990, Karatzas1991}. However, these approaches do not, in general, allow a globally optimal solution.

The theoretic contribution of the paper is to develop a method that gives a globally optimal solution of such mean-variance constrained stochastic  optimal control problems.
More specifically, using the martingale representation theorem,
we reformulate the constraints on the mean and the variance of the agent's payoff function, which are difficult to handle, as two new dynamical systems controlled by new control variables. 
The first new system state represents the agent's future expected payoff with a modified diffusion term. 
The second new system state can be interpreted as the remaining amount of risk the agent can bear. 
The former and latter systems are used to reformulate the constraint on the mean and the variance of the agent's payoff, respectively.
It turns out that the reformulated problem is  risk-sensitive control with a \emph{stochastic target constraint} in the augmented state space of the original system and the new dynamical systems. 
A globally optimal solution to the reformulated problem can be obtained by using the dynamic programming principle in the augmented state space. 
The value function of the problem is computed 
by numerically  solving an associated Hamilton-Jacobi-Bellman equation and is then used to synthesize  optimal compensation and control strategy. 
The proposed solution method allows more general system models for loads and payoff functions for the principal and the agent  than existing dynamic contract methods for the first-best case.
This flexibility and the risk-limiting capability of the proposed contract method make it appropriate for direct load control policy which explicitly treats financial risk in real-time electricity markets.
We also propose an approximate decomposition method for
 the contract design problem:
the problem for $n$ agents can be decomposed into $n$ optimal contract design problems each for a single agent.
This decomposition allows an approximate contract with a provable suboptimality bound.
Due to the decomposition, the computational complexity of the proposed method increases linearly with the number of agents. 
Furthermore, the decomposed contract design problem for an agent is independent of that for another agent. 
Therefore, the contract design procedures for multiple customers can be completely parallelized.


The rest of this paper is organized as follows. 
The problem setting for direct load control and the definition of the risk-limiting dynamic contract are presented in Section \ref{setting}.
We reformulate the constraint on the variance of the agent's payoff, which is called the risk-limiting condition, as a constraint on the compensation provided to the agent by introducing a new control variable in Section \ref{reformulate_compensation}. Using the reformulated constraint, we propose the method for designing a globally optimal risk-limiting dynamic contract 
and discuss its decentralized implementation in Section \ref{design}.
Finally, the performance  of the proposed contract method and its application to direct load control  are demonstrated with data on the electric energy consumption of Austin customers as well as the Electricity Reliability Council of Texas' (ERCOT's) locational marginal price (LMP) data
in Section \ref{application}.

\section{The Setting} \label{setting}

We consider a situation in which the load-serving entity wants to make a contract with $n$  heterogeneous customers to directly control each customer's personal electric load, such as an air conditioner or a water heater. For simplicity, we assume that each customer allows the load-serving entity  control over only one of his or her loads, although the proposed method is also applicable to the case  of multiple loads per customer. 
The load-serving entity's goal is to manage the risk of spending a substantial budget
 in a real-time energy market by controlling the customers' loads in the direct load control program. We consider a finite time horizon contract: let $[0,T]$ be the period in which the contract is effective.

\subsection{Total Power Consumption}

Let $\eta_t^i \in \mathbb{R}$ be the energy consumption (in kWh) up to time $t \in [0,T]$  by customer $i$ and $u_t^i \in \mathbb{R}$ be the power consumption (in kW) by customer $i$'s electric load in the direct load control program. Note that even when $u_t^i = 0$, the total power consumption by customer $i$ is not, in general, zero at time $t$ due to the existence of the customer's other loads and possibly solar or wind generation (which can be considered as negative loads).
If all the customers enter into the contract, the load-serving entity has the authority to determine $u^i:=\{u_t^i\}_{0\leq t \leq T}$ for $i=1, \cdots, n$. 
The number, $n$, of customers is typically in the order of $10^3$--$10^5$.
Let $u := (u^1, \cdots, u^n)$.
The uncertainty in the customers' loads and solar and wind power generation causes the energy consumption process $\{\eta_t\}_{0\leq t \leq T}$ to be stochastic.
To describe the energy consumption process, we use a stochastic differential equation (SDE) model of the form
\begin{equation}\label{demand}
d\eta_t^i = \left (l_i(t) +u_t^i  \right ) dt + \tilde{\sigma}_i(t) dW_t^i,
\end{equation}
where $l_i(t) \in \mathbb{R}$, $0\leq t \leq T$, is the forecast of customer $i$'s loads (in kW) other than those in the direct load control program. 
The effect of the load forecast error is modeled by the diffusion term,  $\tilde{\sigma}_i (t) dW_t^i$, where $W^i := \{W_t^i\}_{0\leq t \leq T}$ is a one-dimensional standard Brownian motion  on a probability space $(\Omega, \mathcal{F}, \mathbb{P})$ and the diffusion coefficient $\tilde{\sigma}_i: [0,T] \to \mathbb{R}$ is a bounded function.
We assume that $W^i$ and $W^j$ are independent for any $i,j \in \{1, \cdots, n\}$ such that $i \neq j$.
The functions $l_i$ and $\tilde{\sigma}_i$ can be estimated from data on the electric energy consumption of customers in Austin as explained in Section \ref{application}.
Furthermore, the validity of the standard Brownian motion in the model for the proposed contract framework is tested using the data in Section \ref{application}.

\subsection{Energy Price and load-serving Entity's Revenue in Real-Time Markets}

Let $p(t) \in \mathbb{R}$, $0 \leq t \leq T$, be the amount of power procured by the load-serving entity in the day-ahead market. We assume that $p(t)$ is given.
The energy price in the real-time market is chosen as the locational marginal price (LMP). Let $\lambda_t$ be the LMP at time $t$. The dynamics of the LMP can be modeled as the following SDE \cite{Deng2001, Kamat2002}:
\begin{equation} \label{mr}
d\lambda_t = r_0 (\nu(t) - \ln \lambda_t)\lambda_t dt + \sigma_0(t) \lambda_t dW_t^0,
\end{equation}
where $W^0 := \{W_t^0\}_{0\leq t \leq T}$ is a one-dimensional standard Brownian motion on $(\Omega, \mathcal{F}, \mathbb{P})$ and the price volatility $\sigma_0: [0, T] \to \mathbb{R}$ is a bounded function.
For simplicity,
we assume that $W^0$ is independent of $W^i$ for $i=1, \cdots, n$, but our contract method can easily be extended to the case in which they are dependent.
This model is suitable to capture the \emph{mean-reverting} behavior of energy prices in the real-time (spot) market: when the energy price is high (resp. low), the supply tends to increase (resp. decrease) and, therefore, causes the price to decrease (resp. increase) \cite{Deng1999}.
Let $w_t:= \ln \lambda_t$, then $w := \{w_t \}_{0\leq t \leq T}$ satisfies
\begin{equation} \label{mr2}
dw_t = r_0 (\nu(t) - w_t) dt + \sigma_0 (t) dW_t^0.
\end{equation}
We estimate $r_0$, $\nu(t)$ and $\sigma_0 (t)$ using the ERCOT LMP data in Section \ref{application}.
In principle, the LMP is not completely exogenous because it is influenced by the power consumption of the customers' loads. In this work, however, we assume that this effect is negligible and, therefore, that the LMP is exogenous.

The load-serving entity's revenue in the real-time market up to time $t$, denoted as $z_t \in \mathbb{R}$, is given by
\begin{equation}\nonumber
\begin{split}
z_t& = \int_0^t \lambda_s \left ( p(s) ds - \sum_{i=1}^n d\eta_s^i \right ).
\end{split}
\end{equation}
Note that we assume that excess power is sold as easily as deficits are procured.
This stochastic integral can be rewritten as the following SDE:
\begin{equation}\label{revenue}
dz_t = \lambda_t \left ( p(t) - \sum_{i=1}^n (l_i (t) + u_t^i) \right ) dt - \sum_{i=1}^n \lambda_t \tilde{\sigma}_i (t) dW_t^i.
\end{equation}
The load-serving entity's revenue is affected by the control $u$ of the customers' loads in the direct load control program.
The set of feasible controls is chosen as
$\mathbb{U}^i := \{ u^i : [0,T] \to \mathcal{U}^i \: | \: u^i \mbox{ progressively measurable with respect to } \mathcal{F}_t^{(i)}  \}$,
where $\mathcal{U}^i$ is a compact set in $\mathbb{R}$ and $\{ \mathcal{F}_t^{(i)} \}_{0\leq t \leq T}$ is the filtration generated by the two dimensional Brownian motion $W^{(i)} := (W^0, W^i)$. 
We also let $\mathbb{U} := \mathbb{U}^1 \times \cdots \times \mathbb{U}^n$.

\subsection{Customers' Loads}

Consider customer $i$'s load in the direct load control program, and
let $x_t^i \in \mathbb{R}$ be the system state at time $t \in [0,T]$. 
If the load is an air conditioner unit, then $x_t^i$ would represent the indoor temperature; if the load is a water heater, it would represent the water temperature.
Then, the system dynamics can be modeled as the following stochastic differential equation:
\begin{equation} \label{indoor}
\begin{split}
dx_t^i =  f_i(x_t^i, u_t^i) dt
\end{split}
\end{equation}
with the initial condition $x_0^i = x^{0i}$ for $i=1, \cdots, n$, where the control $u_t^i$ is a stochastic process. Although our contract method can handle stochastic system models with a diffusion term, we use the model \eqref{indoor}
for simplicity.
Here,  we assume that $f_i : \mathbb{R} \times \mathcal{U}^i \to \mathbb{R}$ is continuous and that $f_i(\bm{x}, \bm{u})$ is differentiable in $\bm{x}$ for any $\bm{u} \in \mathcal{U}^i$. We further assume that there exists a constant $K$ such that for all $(\bm{x}, \bm{u}) \in \mathbb{R} \times \mathcal{U}_i$
\begin{equation} \nonumber
\begin{split}
\left |\frac{\partial f_i (\bm{x}, \bm{u})}{\partial \bm{x}} \right | &\leq K, \\
| f(\bm{x}, \bm{u}) |&\leq K(1 + |\bm{x}| + |\bm{u}|).
\end{split}
\end{equation}
Then, there exists a unique solution $x^i := \{ x_t^i \}_{0\leq t \leq T} \in \mathbb{L}^2(0,T)$ for $i=1, \cdots n$, where 
$\mathbb{L}^2(0,T)$ denotes the space of all real-valued, progressively measurable stochastic processes $\bold{x}$ such that 
$\mathbb{E} \left [ \int_0^T \bold{x}_t^2 dt \right ] < \infty$.
See \cite{Fleming2006} for the proof.

\begin{example} \label{ex_tcl}
Let $x_t^i$ denote the indoor  temperature of customer $i$ at time $t$, and let $\Theta_i(t)$ represent the corresponding outdoor temperature.
Then, the dynamics of the indoor temperature can be described as the following equivalent thermal parameter (ETP) model \cite{Sonderegger1978}:
\begin{equation} \label{indoor_ex}
dx_t^i = [\alpha_i (\Theta_i (t) - x_t^i) - \kappa_i u_t^i] dt
\end{equation}
for $i=1, \cdots, n$. Here, $\alpha_i = R_{1,i}/R_{2,i}$, where $R_{1,i}$ denotes the thermal conductance between the outdoor air and indoor air and $R_{2,i}$ is the thermal conductance between the indoor air and the thermal mass for customer $i$'s room. The positive constant $\kappa_i$ converts an increase in energy (kWh) to a reduction in temperature ($^\circ$C) for customer $i$'s air conditioner.
\end{example}

\subsection{Payoff Functions}

\subsubsection{load-serving entity's payoff}
The load-serving entity's payoff function  is chosen as its profit in the direct load control program. 
Let $C^i \in \mathbb{R}$ be the end-time compensation paid to customer $i$ in the direct load control program and $\mu_i (t)$ be the energy price per unit kWh at time $t$ specified in  customer $i$'s electricity tariff. We assume that $\mu_i : [0,T] \to \mathbb{R}$ is bounded.
The load-serving entity's total payoff in real time, i.e., neglecting the cost of power procured in the day-ahead market, which is its revenue obtained from the customers, is then given by
\begin{equation}  \nonumber
\begin{split}
&\int_0^T dz_t + \sum_{i=1}^n  \int_0^T \mu_i (t) d\eta_t^i - \sum_{i=1}^n C^i \\
&=\sum_{i=1}^n  \int_0^T   \left [ (\mu_i (t) - \lambda_t ) \left ( l_i(t)  + u_t^i \right ) + \lambda_t  p_i(t) \right ] dt+ \sum_{i=1}^n \int_0^T (\mu_i (t) - \lambda_t ) \tilde{\sigma}_i(t) dW_t^i -  \sum_{i=1}^n C^i,
\end{split}
\end{equation}
where $\{p_1(t), \cdots, p_n(t)\}$ is a set satisfying $\sum_{i=1}^n p_i (t) = p(t)$ and
the set of feasible compensation values is chosen as
$\mathbb{C}^i := \{ C^i \in \mathbb{R} \: | \: \mbox{$C^i$ is $\mathcal{F}_T^{(i)}$-measurable} \}$.
Let $\mathbb{C} := \mathbb{C}^1 \times \cdots \times \mathbb{C}^n$.
We define the payoff function of the load-serving entity as
\begin{equation} \label{ppayoff}
\begin{split}
&J^P[C, u] := \sum_{i=1}^n  \left ( \int_0^T r_i^P(t, w_t, x_t^i, u_t^i) dt  + \int_0^T \sigma_i^P (t, w_t) dW_t^i - C^i \right ),
\end{split}
\end{equation}
where $r_i^P: [0,T] \times \mathbb{R}\times \mathbb{R} \times \mathbb{R} \to \mathbb{R}$ and $\sigma_i^P: [0,T] \times \mathbb{R} \to \mathbb{R}$ are such that
\begin{equation} \label{model}
\begin{split}
  r_i^P(t, w_t, x_t^i,u_t^i)  &:= (\mu_i (t) - e^{w_t} )   (u_t^i + l_i (t)) + e^{w_t} p_i(t),\\
\sigma_i^P (t, w_t) &:= (\mu_i (t) - e^{w_t} ) \tilde{\sigma}_i(t)
\end{split}
\end{equation} 
for $i=1, \cdots, n$.  
The superscript `$P$' represents the fact that the load-serving entity plays the role of the principal in the contract.
Note that in the direct load control application $r_i^P$ is independent of $x_t^i$. 
We use the models \eqref{model} for direct load control but the proposed contract design method is applicable to more general models of running payoff and volatility.
For notational simplicity, we will suppress the dependency of the functions on time.

\subsubsection{Customer's payoff}

Each customer's total payoff depends on $(i)$ his or her economic profit and $(ii)$ his or her comfort level.
Customer $i$'s profit can be computed as the compensation received minus the energy costs, i.e.,
\begin{equation}\nonumber
\begin{split}
 &C^i -\int_0^T \mu_i (t) d\eta_t^i \\
 &= -\int_0^T \mu_i (t) (l_i(t) + u_t^i) dt - \int_0^T \mu_i (t)\tilde{\sigma}_i (t) dW_t^i + C^i.
 \end{split}
\end{equation}
 Customer $i$'s payoff function can be represented as
\begin{equation} \label{apayoff}
\begin{split}
J_i^A [C^i, u^i] &:=  \int_0^T r_i^A(x_t^i, u_t^i) dt  + \int_0^T \sigma_i^A(t) dW_t^i +  C^i,
\end{split}
\end{equation}
where $r_i^A : [0,T] \times \mathbb{R} \times \mathbb{R} \to \mathbb{R}$ and $\sigma_i^A : [0,T] \to \mathbb{R}$ are such that
\begin{equation}\label{a_ex}
\begin{split}
r_i^A(t, x_t^i, u_t^i) &:=  -\mu_i (t) (l_i(t) +  u_t^i)  + r_i(x_t^i, u_t^i),\\
\sigma_i^A(t) &:= -\mu_i (t) \tilde{\sigma}_i(t)
\end{split}
\end{equation}
for $i=1, \cdots, n$.  
Here, $r_i(x_t^i, u_t^i)$ represents  customer $i$'s comfort level given the system state $x_t^i$ and control $u_t^i$.
The superscript `$A$' represents the fact that the customer is the agent in the contract.
We assume that there exist constants $B_0$ and $B_1$ such that
$|r_i^A (\bm{x}, \bm{u})| \leq B_0 + B_1 |\bm{x}|$
for all $\bm{x} \in \mathbb{R}$ given any $\bm{u} \in \mathcal{U}^i$.

\begin{example}\label{ex_comfort}
Customer $i$'s discomfort level is zero if the indoor temperature, $x_t^i$, is within a desirable temperature range, $[\underline{\Theta}, \overline{\Theta}]$. 
The discomfort level increases as the indoor temperature increases above $\overline{\Theta}$ or drops below $\underline{\Theta}$. 
If we set the comfort level as the negative value of the discomfort level,
then we can model customer $i$'s comfort level as 
\begin{equation}\label{comfort_tcl}
r_i(x_t^i, u_t^i) = - \omega_i \left [ (x_t^i - \overline{\Theta})_+ + (\underline{\Theta} - x_t^i)_+  \right ],
\end{equation}
where the constant parameter $\omega_i$ represents the customer $i$'s valuation of comfort and $(a)_+ := a$ if $a>0$ and $(a)_+ :=0$ otherwise for any $a \in \mathbb{R}$.
\end{example}

\subsection{Risk-Limiting Dynamic Contracts}

The load-serving entity (principal) offers customer $i$ a contract that specifies the compensation scheme, $C^i$, and its control strategy, $u^i := \{u_t^i \}_{0\leq t \leq T}$ for $i=1, \cdots, n$.
The contract is \emph{dynamic} in the sense that
 the load-serving entity uses the state feedback control strategy written in the contract to dynamically choose the control action each customer must follow. 
Customer $i$ (agent $i$) accepts the contract only if
\begin{enumerate}
\item \label{c1}
(\emph{participation-payoff condition}) 
the mean of customer $i$'s payoff is greater than or equal to some threshold, $b_i \in \mathbb{R}$, i.e.,
\begin{equation}\label{individual}
\mathbb{E} [ J_i^A [C^i, u^i] ]\geq b_i,
\end{equation}
and
\item \label{c2}
 (\emph{risk-limiting condition})  
the variance of customer $i$'s payoff is less than or equal to some threshold, $S_i \in \mathbb{R}$, i.e.,
\begin{equation} \label{risk_limiting}
 \mbox{Var} [ J_i^A[C^i, u^i]] \leq S_i.
 \end{equation}
\end{enumerate}
Note that variance is used as the risk measure of the customer's payoff.
We call $b_i$ and $S_i$ the \emph{participation payoff} and the \emph{risk share} of customer $i$, respectively.

Let $\Lambda := \{(\bold{b}_1, \bold{S}_1), \cdots, (\bold{b}_M, \bold{S}_M) \}$ be a set of given pairs of participation payoffs and risk shares. 
These pairs are designed by the load-serving entity and provided to the customers.
Customer $i$ selects a pair $(b_i, S_i) \in \Lambda$ and the contract determined from this pair.
Once each customer enters into a contract, the load-serving entity directly controls each customer's load following the control strategy specified in the contract. At the end of the contract period, the load-serving entity compensates each customer according to the compensation scheme specified in the contract.

More specifically, once customer $i$ agrees to enter into the contract, the load-serving entity company has the authority to control customer $i$'s load for maximizing the load-serving entity's expected payoff under the constraints \ref{c1}) and \ref{c2}). 
Taking into account the risk of the load-serving entity's payoff being small as well, we formulate the problem of designing such a dynamic contract $(C, u)$ 
as the following constrained \emph{risk-sensitive control} problem:
\begin{subequations} \label{opt}
\begin{align}
\max_{C \in \mathbb{C},u \in \mathbb{U}} \quad & - \frac{1}{\theta} \log \mathbb{E}\left [\exp(-\theta J^P[C, u])  \right] \label{obj} \\
\mbox{subject to} \quad &dw_t = r_0 (\nu(t) - w_t) dt + \sigma_0 (t) dW_t^0 \\
& dx_t^i = f_i(x_t^i, u_t^i) dt \label{temp}\\
& \mathbb{E} [J_i^A[C^i, u^i]] \geq b_i \label{ind}\\
& \mbox{Var} [ J_i^A[C^i, u^i]] \leq S_i,\label{ris} 
\end{align}
\end{subequations}
where $\theta \in \mathbb{R} \setminus \{0\}$ is a constant, called the coefficient of load-serving entity's \emph{risk-aversion}.
When $\theta$ is positive,
the risk-sensitive objective function penalizes the risk of the load-serving entity's payoff being small because $-\exp(-  \theta J^P)$ is concave increasing in $J^P$. Therefore, the load-serving entity can make a risk-averse decision by solving \eqref{opt}. 
If $\theta < 0$, the load-serving entity is risk-seeking. For intuition, note that the risk-sensitive objective function is well approximated by a weighted sum of the mean and the variance of the payoff when $|\theta |$ is small because the Taylor expansion of the risk-sensitive objective function is given by
\begin{equation}\label{interp_rs}
-\frac{1}{\theta} \log \mathbb{E} \left[\exp ( -\theta J^P)  \right ] = \mathbb{E} [J^P ] - \frac{\theta}{2} \mbox{Var} [ J^P] + O(\theta^2)
\end{equation}
as $\theta \to 0$. Note that the variance of the payoff is penalized when $\theta > 0$.

The solution, $(C^{\mbox{\tiny OPT}}, u^{\mbox{\tiny OPT}})$, to this problem is said to be the optimal \emph{risk-limiting dynamic contract} for direct load control.
This problem of risk-limiting dynamic contract design is a mean-variance constrained-stochastic optimal control problem that is not directly solvable via dynamic programming. 
In the following sections, we carefully characterize the necessary and sufficient conditions for the constraints on the mean and the variance of the customers' payoff functions. The characterizations allow us to show that the contract design problem can be reformulated as a risk-sensitive control problem with a stochastic target constraint that can be solved by dynamic programming.

 The information availabilities to the load-serving entity and the customers in the direct load control program are symmetric, as opposed to the case of indirect load control, in which the load-serving entity has limited observation capability (e.g., \cite{Yang2014}). More specifically, in the proposed framework, the load-serving entity can monitor the control and state of the customers' loads in the direct load control program as well as the energy price in the real-time market.
 Furthermore, the load-serving entity has all the parameters and functions needed to design an optimal contract. That is, it has the information of $p$, $l_i$, $\tilde{\sigma}_i$, $\mu_i$, $r_0$, $\sigma_0$, $\nu$, which can be estimated from data as shown in Section \ref{application}, and the customers' comfort functions $r_i$ and load models for $i=1, \cdots, n$. In practice, the comfort functions and load models can be identified using a training period.
 Each customer, in principle, can have the same information. 
However,
customer $i$ needs to know only $p$, $l_i$, $\tilde{\sigma}_i$, $\mu_i$, $r_0$, $\sigma_0$, $\nu$, his or her own comfort function and load model.
The proposed framework assumes that each customer can monitor the control and state of his or her load in the direct load control program and the energy price in the real-time market.
 Another important feature of the proposed contract method is that the interactions between the load-serving entity and its customers can be decoupled from each other because one customer's load does not affect those of other customers and the participation payoff and risk-limiting conditions are personalized. 
This feature will allow us to decentralize the control of loads
 as shown in Section \ref{multi}.

We assume that the contract period $[0,T]$ is a time interval within 24 hours, but the proposed method can handle arbitrary finite time horizons. Therefore, the customers and the load-serving entity, in principle, can renew the contracts every day. 
However, it may not be convenient for each customer to choose a contract or, equivalently, a participation payoff and a risk share pair every day. 
This issue can be resolved by automatically choosing the contract for the current day as that for the previous day unless the customer explicitly wants to change it.
Daily contracts have a practical advantage: the day-ahead forecasts of the LMP model parameters and demand uncertainty (and outdoor temperature in the case of air conditioners) can be incorporated into the contracts. Therefore, the contracts can be designed using accurate models.
We also assume that each customer does not strategically control other loads to modify the forecasted $\tilde{\sigma}_i$ by the load-serving entity. This assumption can be justified in two ways. 
First, the load-serving entity can make a contract to control multiple loads of a customer so that the customer has little flexibility to change $\tilde{\sigma}_i$. 
Second, even if the customer strategically affects $\tilde{\sigma}_i$ in one day, the customer's gain in the next day is marginal because the contract is renewed with a new estimate $\tilde{\sigma}_i$ that incorporates any strategic behavior.  
Formally, this problem can be formulated as a Stackelberg differential game, in which the the load-serving entity chooses the estimates of $l_i$ and $\tilde{\sigma}_i$ for the contract period $[0,T]$ assuming that the customer has no incentive to deviate from $\tilde{\sigma}_i$ in the contract period.
A similar problem is considered in our previous work \cite{Yang2014}.
This Stackelberg differential game problem is out of the scope of this paper and will be addressed in our future work on risk-limiting dynamic contracts for indirect load control.

\section{Risk-Limiting Compensation} \label{reformulate_compensation}

The risk-limiting condition \eqref{risk_limiting} is an inequality constraint on the variance of each agent's payoff. This constraint hinders us from using the dynamic programming principle to solve the contract design problem \eqref{opt}. 
In this section, we characterize a condition on the end-time compensation, which is equivalent to the risk-limiting condition.
It turns out that the new equivalent condition allows us to formulate the contract design problem as a risk-sensitive control problem that can be solved by dynamic programming.

We begin by defining the following set of stochastic processes: let $\Gamma^i$ be the set of processes $\xi^i:= \{\xi_t^i\}_{0\leq t \leq T}$, $\xi_t^i = (\xi_t^{i,1}, \xi_t^{i,2}) \in \mathbb{R}^{1 \times 2}$ such that
\begin{enumerate}[(i)]
\item $\xi_t^i$ is $\mathcal{F}_t^{{(i)}}$-progressively measurable; 
\item $\mathbb{E} \left [ \int_0^T \| \xi_t^i \|^2 dt  \right ] < \infty$
\end{enumerate}
for $i=1, \cdots, n$.
We also let $\Gamma := \Gamma^1\times \cdots \times \Gamma^n$. 
In the next lemma, we show that there exists a unique process in this set such that its integral over the Brownian motion $W^{(i)}:= (W^0, W^i)$ corresponds to the difference between the agent's payoff and its mean value.

\begin{lemma} \label{lem1}
Fix $i \in \{1, \cdots, n\}$.
Given $C^i \in \mathbb{C}^i$ and $u^i \in \mathbb{U}^i$, there exists a unique (up to a set of measure zero) stochastic process $\xi^i = \{\xi_t^i \}_{0\leq t \leq T} \in \Gamma^i$ such that 
\begin{equation}\label{lem1_eq}
\begin{split}
&J_i^A [ C^i, u^i] - \mathbb{E}[ J_i^A [ C^i, u^i] ] =   \int_0^T \xi_t^i  dW_t^{(i)}.
\end{split}
\end{equation}
\end{lemma}
\begin{IEEEproof} 
Fix $C^i \in \mathbb{C}^i$ and $u^i \in \mathbb{U}^i$.
We introduce a new process
\begin{equation}\nonumber
q_t^i := \mathbb{E} \left [ \left. \int_t^T r_i^A(x_s^i, u_s^i) ds  + C^i \: \right | \: \mathcal{F}_t^{(i)}  \right ].
\end{equation}
Here, the expectation is conditioned over the filtration $\{\mathcal{F}_t^{(i)} \}_{0\leq t \leq T}$ generated by the Brownian motion $W^{(i)} = (W^0, W^i)$. We notice that the process
\begin{equation} \nonumber
\begin{split}
&q_t^i + \int_0^t r_i^A(x_s^i, u_s^i) ds = \mathbb{E} \left [  \left. \int_0^T r_i^A(x_s^i, u_s^i) ds + C^i  \: \right | \: \mathcal{F}_t^{(i)} \right ]
\end{split}
\end{equation}
is martingale. 
Recall that there exist constants $B_0$ and $B_1$ such that
\begin{equation}\nonumber
\left | r_i^A (\bm{x}, \bm{u})  \right | \leq B_0 + B_1 | \bm{x} | 
\end{equation}
for all $\bm{x} \in \mathbb{R}$ and $\bm{u} \in \mathcal{U}^i$. Therefore, we have
\begin{equation}\nonumber
\left ( \int_0^t r_i^A(x_s^i, u_s^i) ds \right )^2 \leq \int_0^t (B_0 + B_1 |x_s^i|)^2 ds, 
\end{equation}
which implies that for $t \in [0,T]$
\begin{equation} \label{ineq1}
\mathbb{E} \left [ \left ( \int_0^t r_i^A(x_s^i, u_s^i) ds \right )^2  \right ] < \infty
\end{equation}
because $x^i \in \mathbb{L}^2(0,T)$.
From the definition of $q_t^i$, we deduce that
\begin{equation} \label{ineq2}
\mathbb{E} \left [   (q_t^i)^2  \right ] < \infty.
\end{equation} 
Due to the inequalities \eqref{ineq1} and \eqref{ineq2}, we obtain
\begin{equation}\nonumber
\mathbb{E} \left [  \left ( q_t^i + \int_0^t r_i^A (x_s^i, u_s^i)ds \right )^2 \right ]< \infty.
\end{equation}
The Martingale representation theorem (e.g., \cite{Karatzas1998, Korn2001}) suggests that there exists a unique (up to set of measure zero) process $\bar{\xi}^i = \{\bar{\xi}_t^i \}_{0\leq t \leq T} \in \Gamma^i$ such that
\begin{equation} \nonumber
\begin{split}
&q_T^i  + \int_0^T r_i^A(x_t^i, u_t^i) dt  = q_0^i + \int_0^T \bar{\xi}_t^i  dW_t^{(i)}.
\end{split}
\end{equation}
We also note that
\begin{equation}\nonumber
\begin{split}
q_T^i &=  C^i,\\
q_0^i &= \mathbb{E} [ J_i^A [C^i, u^i]].
\end{split}
\end{equation}
Therefore, we obtain
\begin{equation}\nonumber
\begin{split}
J_i^A [ C^i, u^i] &= \mathbb{E}[ J_i^A [ C^i, u^i] ]+ \int_0^T  \sigma_i^A(t) dW_t^i
+ \int_0^T \bar{\xi}_t^{i}  dW_t^{(i)}.
\end{split}
\end{equation}
Set $\xi_t^{1,i} := \bar{\xi}_t^{1,i}$ and
$\xi_t^{2,i} := \bar{\xi}_t^{2,i} + \sigma_i^A(t)$, then $\xi^i$ is in $\Gamma^i$ and satisfies \eqref{lem1_eq}.
\end{IEEEproof}
This lemma represents the agent's payoff as the sum of its mean value 
and the It\^{o} integral of the new process $\xi^i$ along the Brownian motion $W^{(i)}$.
The following theorem suggests that this representation allows reformulation of the risk-limiting condition \eqref{risk_limiting}.
\begin{theorem}\label{thm1}
Fix $i \in \{1, \cdots, n\}$ and $u^i \in \mathbb{U}^i$.
The risk-limiting condition holds, i.e.,
\begin{equation} \label{rs_thm}
\mbox{\emph{Var}} [ J_i^A[C^i, u^i] ] \leq S_i
\end{equation}
 if and only if there exists a unique (up to set of measure zero) $\gamma^i \in \Gamma^i$ such that
\begin{equation} \label{condition1}
\begin{split}
C^i &= \mathbb{E} [ J_i^A[C^i, u^i]] - \int_0^T r_i^A (x_t^i, u_t^i) dt - \int_0^T \sigma_i^A(t) dW_t^i + \int_0^T \gamma_t^i  dW_t^{(i)}
\end{split}
\end{equation}
and
\begin{equation} \label{condition2}
\mathbb{E} \left [  \int_0^T \| \gamma_t^i \|^2dt \right] \leq S_i.
\end{equation}
\end{theorem}
\begin{IEEEproof}
Suppose that there exists $\gamma^i \in \Gamma^i$ such that \eqref{condition1} and \eqref{condition2} hold. Then,
\begin{equation}\nonumber
J_i^A[C^i, u^i] - \mathbb{E} [ J_i^A[C^i, u^i] ] = \int_0^T \gamma_t^i dW_t^{(i)}.
\end{equation}
Due to the It\^{o}'s isometry, we have
\begin{equation}\nonumber
\mbox{Var} [ J_i^A[C^i, u^i] ] =\mathbb{E} \left [  \int_0^T \|\gamma_t^i \|^2 dt \right].
\end{equation}
Combining this equality and \eqref{condition2}, we obtain 
 the risk-limiting condition \eqref{rs_thm}. 

Suppose now that the risk-limiting condition \eqref{rs_thm} holds.
 Lemma \ref{lem1} suggests that there exists $\xi^i \in \Gamma^i$ such that
\begin{equation}\nonumber
\begin{split}
&J_i^A [ C^i, u^i] - \mathbb{E}[ J_i^A [ C^i, u^i] ] = \int_0^T \xi_t^i dW_t^{(i)}.
\end{split}
\end{equation}
The variance of agent $i$'s payoff is given by
\begin{equation}\nonumber
\mbox{Var} [ J_i^A[C^i, u^i]] = \mathbb{E} \left [  \int_0^T \| \xi_t^i \|^2 dt \right ].
\end{equation}
Due to the risk-limiting condition \eqref{rs_thm}, we have
\begin{equation} \nonumber
\mathbb{E} \left [ \int_0^T \| \xi_t^i \|^2 dt \right ] \leq S_i.
\end{equation}
Therefore, $\xi^i \in \Gamma^i$ satisfies both \eqref{condition1} and \eqref{condition2}.
Such a process must be unique up to a set of measure zero due to Lemma \ref{lem1}.
\end{IEEEproof}

The next corollary suggests a way to construct the end-time compensation given $u \in \mathbb{U}$ and $\gamma \in \Gamma$.

\begin{cor} \label{cor1}
Fix $u \in \mathbb{U}$ and $\gamma \in \Gamma$ such that 
\begin{equation} \label{equiv_condition}
\mathbb{E} \left[\int_0^T \| \gamma_t^i \|^2 dt  \right ] \leq S_i
\end{equation}
for $i=1, \cdots, n$.
The risk-limiting condition \eqref{risk_limiting} holds if and only if the end-time compensation, $C \in \mathbb{C}$, satisfies
\begin{equation} \nonumber 
\begin{split}
C^i &= \mathbb{E} [ J_i^A[C^i, u^i]] - \int_0^T r_i^A (x_t^i, u_t^i) dt 
  - \int_0^T \sigma_i^A(t) dW_t^i+ \int_0^T \gamma_t^i dW_t^{(i)}
\end{split}
\end{equation}
for $i=1, \cdots, n$. 
\end{cor}
Theorem \ref{thm1} and Corollary \ref{cor1} imply that,  given $u^i \in \mathbb{U}^i$, determining $C^i \in \mathbb{C}^i$ is equivalent to choosing $\gamma^i \in \Gamma^i$ such that it satisfies \eqref{equiv_condition}.
In the next section, we consider $\gamma^i \in \Gamma^i$ as a two-dimensional decision variable and then construct the end-time compensation using an optimal $\gamma^i$. 
We also show that the mean of agent $i$'s payoff is given by the participation payoff $b_i$ if an optimal contract is chosen.
However, even if we reformulate the contract design problem \eqref{opt} as a stochastic optimal control in which the decision variables are  $u^i$ and $\gamma^i$, $i=1, \cdots, n$, the integral constraint \eqref{equiv_condition} prohibits us from using dynamic programming to solve the reformulated problem.
We resolve this issue  in the following section by  introducing new state and control variables.

\section{Risk-Limiting Dynamic Contract Design} \label{design}

We now propose the solution method for risk-limiting dynamic contract design problem \eqref{opt} given $(b_i, S_i) \in \Lambda$ for $i=1, \cdots, n$.
Recall that the principal's objective is to maximize its risk-sensitive payoff with a guarantee that the participation payoff and risk-limiting conditions for all the agents are satisfied. The risk of each agent's payoff being too small is limited by the variance constraint \eqref{ris}, which is the risk-limiting condition.
On the other hand, small values of the principal's payoff are penalized by maximizing the risk-sensitive objective function \eqref{obj} when the risk-aversion coefficient $\theta$ is positive.
In other words, if the principal is risk-averse, she transfers her financial risk to the agents as long as the agents' risk-limiting conditions are respected.

The contract design problem \eqref{opt} is a constrained stochastic optimal control problem that cannot be directly solved by dynamic programming \cite{Bellman1956}.
The stochastic maximum principle approach may be used to handle the constraints \cite{Peng1990, Cadenillas1995}. However, it is not capable of finding a globally optimal solution unless the problem is concave  in both $C$ and $u$. 
To obtain a globally optimal solution, we reformulate the problem as a risk-sensitive control problem that can be solved by dynamic programming. 
The key idea is to introduce two new state variables. 
The first state variable's value at the terminal time allows us to construct the end-time compensation using Theorem \ref{thm1}.
The second variable is used to reformulate the constraint \eqref{condition2} on the new decision variable, $\gamma$, as an SDE.
However, the dynamic programming approach, in general, has an inherent scalability issue: the computational complexity exponentially increases as the system dimension increases (e.g., \cite{Bertsekas2005}). 
We overcome this scalability issue by proposing an approximate decomposition of the contract design problem for all agents into $n$ lower-dimensional contract design problems, each for a single agent, where $n$ is the number of agents.

\subsection{Reformulation} 

We show that the solution of the contract design problem \eqref{opt} can be obtained by solving the following risk-sensitive control problem:
\begin{subequations} \label{simple}
\begin{align}
\max_{u \in \mathbb{U}, \gamma \in \Gamma, \zeta \in \Gamma} \quad &-\frac{1}{\theta} \log \mathbb{E}\left [\exp \left (-\theta \bar{J}^P [u, \gamma, \zeta]  \right )  \right] \\
\mbox{subject to} \quad  &dw_t = r_0 (\nu(t) - w_t) dt + \sigma_0 (t) dW_t^0\\
& dx_t^i = f_i(x_t^i, u_t^i) dt \label{temp_s} \\
& dv_t^i = -r_i^A(x_t^i, u_t^i) dt + \gamma_t^{i,1} dW_t^0 + (\gamma_t^{i,2}  - \sigma_i^A(t)) dW_t^i \label{v_s}  \\ 
&v_0^i = b_i  \label{init_s} \\ 
& dy_t^i = - \| \gamma_t^i \| ^2 dt +  \zeta_t^i dW_t^{(i)} \label{sde2}\\
& y_0^i = S_i \label{init2}\\
& y_T^i \geq 0 \quad \mbox{a.s.}, \label{ter2}
\end{align}
\end{subequations}
where 
$\bar{J}^P$ is the reformulated principal's payoff given by
\begin{equation} \nonumber
\begin{split}
&\bar{J}^P [u, \gamma, \zeta] := \sum_{i=1}^n  \left (\int_0^T r_i^P(w_t, x_t^i, u_t^i) + \int_0^T \sigma_i^P(w_t) dW_t  - v_T^i \right ).
\end{split}
\end{equation}
Note that we now view $\gamma \in \Gamma$ as a decision variable instead of $C \in \mathbb{C}$. This is feasible due to Theorem \ref{thm1} and Corollary \ref{cor1}.
The new state $y^i$ and new decision variable $\zeta^i$ handle the integral constraint \eqref{condition2} in Theorem \ref{thm1}.
Intuitively speaking, 
the first new state variable $v_t^i$ represents agent $i$'s expected future payoff  with a modified diffusion term.
The second new state $y_t^i$ can be interpreted as the remaining amount of risk that agent $i$ can bear from time $t$.
Having these interpretations, we can show that
the terminal value of the first new state variable can be used to construct an optimal end-time compensation
and that of the second new  state variable must be greater than or equal to zero to satisfy agent $i$'s risk-limiting condition.
These two claims are shown in Theorem \ref{optimal}.

Another important observation is that the problem is now defined in the augmented state space of $w_t$, $x_t := (x_t^1, \cdots, x_t^n) \in \mathbb{R}^n$, $v_t := (v_t^1, \cdots, v_t^n) \in \mathbb{R}^n$ and $y_t := (y_t^1, \cdots, y_t^n) \in \mathbb{R}^n$. Therefore, the total system dimension is $3n+1$.
This reformulated problem 
can be decentralized into $n$  three dimensional risk-sensitive control problems, as shown in Section \ref{multi}.

\begin{theorem} \label{optimal}
Let $(u^{\mbox{\tiny OPT}}, \gamma^{\mbox{\tiny OPT}}, \zeta^{\mbox{\tiny OPT}})$ be the solution to \eqref{simple}. 
We also let $x^{\mbox{\tiny OPT}}$, $v^{\mbox{\tiny OPT}}$ and $y^{\mbox{\tiny OPT}}$ denote the processes driven by \eqref{temp_s}, \eqref{v_s} and \eqref{sde2} with $(u^{\mbox{\tiny OPT}}, \gamma^{\mbox{\tiny OPT}}, \zeta^{\mbox{\tiny OPT}})$, respectively.
Define
\begin{equation} \label{comp_opt}
\begin{split}
C^{\mbox{\tiny OPT},i} &:= v_T^{\mbox{\tiny OPT},i}
 \end{split}
\end{equation}
for $i =1, \cdots, n$.
Then, $(C^{\mbox{\tiny OPT}}, u^{\mbox{\tiny OPT}})$ is an optimal risk-limiting dynamic contract, i.e., it solves \eqref{opt}.
\end{theorem}

\begin{IEEEproof}
We first observe that 
\begin{equation} \label{expression}
\begin{split}
J_i^A[C^{\mbox{\tiny OPT},i}, u^{\mbox{\tiny OPT},i}] &= \int_0^T r_i^A(x_t^{\mbox{\tiny OPT},i}, u_t^{\mbox{\tiny OPT},i}) dt + \int_0^T \sigma_i^A(t) dW_t^i  + v_T^{*i} \\
&= b_i + \int_0^T \gamma_t^{\mbox{\tiny OPT},i}  dW_t^{(i)}
\end{split}
\end{equation}
due to the SDE \eqref{v_s} with the initial condition \eqref{init_s}.
Therefore, we have
\begin{equation} \label{binding}
\mathbb{E} [J_i^A[C^{\mbox{\tiny OPT},i}, u^{\mbox{\tiny OPT},i}]] = b_i,
\end{equation}
which implies that the participation payoff condition \eqref{ind} holds. Furthermore, the variance of agent $i$'s payoff with $(C^{\mbox{\tiny OPT},i}, u^{\mbox{\tiny OPT},i})$ is given by
\begin{equation} \nonumber
\mbox{Var} [J_i^A[C^{\mbox{\tiny OPT},i}, u^{\mbox{\tiny OPT},i}]] = \mathbb{E} \left [  \int_0^T \| \gamma_t^{\mbox{\tiny OPT},i} \|^2 dt \right ]
\end{equation}
due to the It\^{o}'s isometry.
We also notice that
\begin{equation} \nonumber
y_T^{\mbox{\tiny OPT},i} = S_i - \int_0^T \| \gamma_t^{\mbox{\tiny OPT},i} \|^2 dt + \int_0^T \zeta_t^{\mbox{\tiny OPT},i} dW_t^{(i)},
\end{equation}
which suggests that
\begin{equation} \nonumber
\mbox{Var} [J_i^A[C^{\mbox{\tiny OPT},i}, u^{\mbox{\tiny OPT},i}]] = \mathbb{E} \left [ S_i - y_T^{\mbox{\tiny OPT},i} + \int_0^T \zeta_t^{\mbox{\tiny OPT},i} dW_t^{(i)} \right].
\end{equation}
Hence, if $y_T^{\mbox{\tiny OPT},i} \geq 0$ a.s., the risk-limiting condition \eqref{ris} also holds.
Therefore, $(C^{\mbox{\tiny OPT}}, u^{\mbox{\tiny OPT}})$ is a feasible dynamic contract.

Suppose that $(C^{\mbox{\tiny OPT}}, u^{\mbox{\tiny OPT}})$ is not a solution of \eqref{opt} and select a solution, $(\hat{C}, \hat{u})$, of \eqref{opt}.
Theorem \ref{thm1} suggests that there exists a unique (up to set of measure zero) $\hat{\gamma} \in \Gamma$ such that 
\begin{equation} \nonumber
\begin{split}
\hat{C}^i &= \mathbb{E} [ J_i^A [ \hat{C}^i, \hat{u}^i]] - \int_0^T r_i^A (\hat{x}_t^i, \hat{u}_t^i) dt -\int_0^T \sigma_i^A(t) dW_t^i + \int_0^T \hat{\gamma}_t^i  dW_t^{(i)}
\end{split}
\end{equation}
and 
\begin{equation} \label{ris_condition}
\mathbb{E} \left [ \int_0^T \| \hat{\gamma}_t^i \|^2 dt \right ] \leq S_i.
\end{equation}
We claim that $(\hat{u}, \hat{\gamma})$ satisfies all the constraints of \eqref{simple} with the 
 stochastic process $\hat{v} := \{ \hat{v}_t \}_{0\leq t\leq T}$ defined as
\begin{equation}\nonumber
\begin{split}
\hat{v}_t^i &:= \hat{v}_0^i - \int_0^t r_i^A (\hat{x}_s^i, \hat{u}_s^i) ds + \int_0^t (\hat{\gamma}_s^i - \sigma_i^A(s) ) dW_s^i
\end{split}
\end{equation}
for some initial value $\hat{v}_0^i \in \mathbb{R}$ such that
\begin{equation} \label{eq2}
\hat{v}_T^i = \hat{C}^{i}.
\end{equation}
It is clear that the process $\hat{v}^i$ satisfies the SDE \eqref{v_s} for $i=1, \cdots, n$ by definition.
Additionally, note that
\begin{equation} \nonumber
\mathbb{E} [J_i^A[\hat{C}^i, \hat{u}^i] ] = \hat{v}_0^i.
\end{equation}
Suppose that $\hat{v}$ does not satisfy the initial condition \eqref{init_s}. 
We first assume that there exists $j \in \{1, \cdots, n\}$ such that $\hat{v}_0^j > b_j$. 
Define a new end-time compensation $C'$ as 
\begin{equation}\nonumber
C'^{i} := \left \{
\begin{array}{ll}
\hat{C}^j  - (\hat{v}_0^j - b_j) & \mbox{ if $i = j$}\\
\hat{C}^i & \mbox{ otherwise}.
\end{array}
\right.
\end{equation}
We then have
\begin{equation}\nonumber
\begin{split}
\mathbb{E} [J_j^A [C'^j, \hat{u}^j]] - C'^j &= \mathbb{E} [J_j^A [\hat{C}^j, \hat{u}^j]] - \hat{C}^j\\
&= \hat{v}_0^j - \hat{C}^j,
\end{split}
\end{equation}
which implies that
\begin{equation}\nonumber
\mathbb{E} [ J_j^A [ C'^{j}, \hat{u}^j]] = b_j.
\end{equation}
Therefore, $(C', \hat{u})$ satisfies the participation payoff condition \eqref{ind} for all $i=1, \cdots, n$. 
The risk-limiting condition also holds with $(C', \hat{u})$ because the difference between $C'$ and $\hat{C}$ is deterministic.
On the other hand, we notice that 
\begin{equation}\nonumber
J^P [ C', \hat{u}] > J^P [ \hat{C}, \hat{u}]
\end{equation}
because $C'^j < \hat{C}^j$ and $C'^i = \hat{C}^i$ for $i\neq j$. 
The contract 
$(C', \hat{u})$ satisfies all the constraints of \eqref{opt} 
and is strictly better than $(\hat{C}, \hat{u})$. 
This is a contradiction because $(\hat{C}, \hat{u})$ solves the contract design problem \eqref{opt}.
Therefore, $\hat{v}_0^i$ must be equal to the participation payoff $b_i$ for $i=1, \cdots, n$. Hence, $(\hat{u}, \hat{\gamma})$ satisfies all the constraints of  \eqref{simple} with the processes $\hat{y}$ and $\hat{x}$, where $\hat{x}$ solves \eqref{temp_s} with the control $\hat{u}$.

We define a stochastic process $\tilde{y}^i := \{\tilde{y}_t^i\}_{0\leq t \leq T}$ as
\begin{equation} \nonumber
\tilde{y}_t^i := \mathbb{E} \left [ \left. \int_t^T \|  \hat{\gamma}_t^i\| ^2 dt \: \right | \: \mathcal{F}_t^{(i)} \right ].
\end{equation}
Note that 
\begin{equation}\nonumber
\tilde{y}_t^i + \int_0^t \|\hat{\gamma}_s^i \|^2 ds = \mathbb{E} \left [ \left. \int_0^T \|\hat{\gamma}_t^i \|^2 dt \: \right | \: \mathcal{F}_t^{(i)} \right ]
\end{equation}
is martingale. Furthermore, $\mathbb{E} \left [ \left ( \tilde{y}_t^i + \int_0^t \| \hat{\gamma}_s^i \|^2  ds\right )^2 \right ] < \infty$ because $\hat{\gamma}^i \in \Gamma^i$.
Therefore, the martingale representation theorem suggests that 
there exists a unique (up to a set of measure zero) $\hat{\zeta}^i \in \Gamma^i$ such that 
\begin{equation}\nonumber
\tilde{y}_t^i  + \int_0^t \|\hat{\gamma}_t^i \|^2 dt= \tilde{y}_0^i + \int_0^t \hat{\zeta}_t^i dW_t^{(i)}.
\end{equation}
Therefore, the process $\tilde{y}^i$ solves \eqref{sde2} with $(\hat{\gamma}, \hat{\zeta})$.
Due to \eqref{ris_condition}, we also have 
\begin{equation} \nonumber
\tilde{y}_0^i \leq S_i
\end{equation}
 and hence the constraint \eqref{init2} is satisfied.
We define another stochastic process $\hat{y}^i := \{\hat{y}_t^i\}_{0\leq t \leq T}$ as
\begin{equation} \nonumber
\hat{y}_t^i := \tilde{y}_t^i + S_i - \tilde{y}_0^i
\end{equation}
for $t \in [0,T]$.
Therefore, we have
$\hat{y}_0^i = S_i$  and $\hat{y}_T^i \geq 0$
because $\tilde{y}_T^i = 0$ by definition. Hence, $\hat{y}^i$ satisfies the constraints \eqref{init2} and \eqref{ter2}.

Since $(\hat{C}, \hat{u})$ solves \eqref{opt}, while $(C^{\mbox{\tiny OPT}}, u^{\mbox{\tiny OPT}})$ does not, the following inequality holds:
\begin{equation}\nonumber
\mathbb{E}  [ -\exp ( -\theta J^P[\hat{C}, \hat{u}] )  ] > \mathbb{E}  [ -\exp( -\theta J^P[C^{\mbox{\tiny OPT}}, u^{\mbox{\tiny OPT}}] )  ].
\end{equation}
This inequality can be rewritten as
\begin{equation}\nonumber
\begin{split}
&\mathbb{E}  \left [- \exp \left (-\theta \sum_{i=1}^n \left ( \int_0^T r_i^P(w_t, \hat{x}_t^i, \hat{u}_t^i) dt -  \hat{v}_T^i \right ) \right )\right ] >\\
&\mathbb{E}  \left [- \exp \left (-\theta \sum_{i=1}^n \left ( \int_0^T r_i^P(w_t, x_t^{\mbox{\tiny OPT},i}, u_t^{\mbox{\tiny OPT},i}) dt -  v_T^{\mbox{\tiny OPT},i} \right ) \right ) \right ]
\end{split}
\end{equation}
due to \eqref{eq2} and \eqref{comp_opt}.
This is contradictory to the fact that $(u^{\mbox{\tiny OPT}}, \gamma^{\mbox{\tiny OPT}}, \zeta^{\mbox{\tiny OPT}})$ is a solution to \eqref{simple}.
Therefore, $(C^{\mbox{\tiny OPT}}, u^{\mbox{\tiny OPT}})$ should solve \eqref{opt} and hence an optimal risk-limiting dynamic contract.
\end{IEEEproof}
We observe that the agents' expected payoff must be equal to their participation payoffs from \eqref{binding}, i.e., the inequalities \eqref{ind} for the participation payoff are always binding at an optimal contract. 
Intuitively speaking, if the agent's expected payoff is strictly greater than his or her participation payoff,
the principal has an incentive to decrease the end-time compensation for the agent.
The equality \eqref{expression} also suggests that each agent's payoff can be completely characterized by his or her participation payoff and the new control variable $\gamma$ if an optimal contract is executed.
\begin{cor}\label{agent_payoff}
Agent $i$'s payoff with an optimal risk-limiting dynamic contract $(C^{\mbox{\tiny OPT}}, u^{\mbox{\tiny OPT}})$ is given by
\begin{equation} \nonumber
J_i^A[C^{\mbox{\tiny OPT},i}, u^{\mbox{\tiny OPT},i}] = b_i + \int_0^T \gamma_t^{\mbox{\tiny OPT},i} dW_t^{(i)}
\end{equation}
for $i = 1, \cdots, n$. Therefore, 
\begin{equation} \nonumber
\mathbb{E} \left [ J_i^A[C^{\mbox{\tiny OPT},i}, u^{\mbox{\tiny OPT},i}] \right ] = b_i.
\end{equation}
\end{cor}

\subsection{Decoupled Contract Design and Decentralized Control} \label{multi}

 We propose an approximate decomposition of the contract design problem \eqref{opt} into $n$ low dimensional problems using the fact that the system dynamics \eqref{temp}, the participation payoff condition \eqref{ind} and the risk-limiting condition \eqref{ris} for one agent are decoupled from those for other agents and that $W^1, \cdots, W^n$ are mutually independent. The approximate solution obtained using this decomposition has a guaranteed suboptimality bound.
This decomposition enables the direct load control program with the proposed dynamic contracts to handle a large population of agents without scalability issues.

More specifically, for each $i \in \{1,\cdots, n\}$, the  approximate risk-limiting dynamic contract for agent $i$ can be obtained by solving the following  risk-sensitive control problem:
\begin{equation} \label{dec}
\begin{split}
\max_{\substack{u^i \in \mathbb{U}^i,\\ \gamma^i \in \Gamma^i, \zeta^i \in \Gamma^i}} \quad &-\frac{1}{\theta}\log \mathbb{E}\left [ \exp \left( -\theta \bar{J}_i^P [ u^i, \gamma^i, \zeta^i] \right )\right]  \\
\mbox{subject to} \quad & 
dw_t = r_0 (\nu(t) - w_t) dt + \sigma_0 (t) dW_t^0 \\
&dx_t^i = f_i(x_t^i, u_t^i) dt \\
& dy_t^i = -\| \gamma_t^i \| ^2 dt + \zeta_t^i dW_t^{(i)} \\
& y_0^i = S_i \\
&y_T^i \geq 0 \quad \mbox{a.s.},
\end{split}
\end{equation}
where 
\begin{equation} \nonumber
\begin{split}
\bar{J}_i^P [u^i, \gamma^i, \zeta^i] &:=  -b_i +\int_0^T (r_i^P(w_t, x_t^i, u_t^i ) + r_i^A(x_t^i, u_t^i)) dt+ \int_0^T (\sigma_i^P (w_t) + \sigma_i^A(t) ) dW_t^i
- \int_0^T \gamma_t^{i} dW_t^{(i)}.
 \end{split}
\end{equation}
Note that the system for $v^i$ is absorbed into the modified payoff function $\bar{J}_i^P$.
 It is clear that this decomposition is exact when $\sigma_0 \equiv 0$ due to the mutual independence of $\{W^1, \cdots, W^n\}$.
In addition, the following proposition suggests that
the approximate contract obtained using the proposed decomposition has a provable suboptimality bound, which can be computed {\it a posteriori}.
\begin{proposition}
Let $(u^{\mbox{\tiny OPT}}, \gamma^{\mbox{\tiny OPT}}, \zeta^{\mbox{\tiny OPT}})$ and $(u^{*i}, \gamma^{*i}, \zeta^{*i})$ be the solutions to \eqref{opt} and \eqref{dec}, respectively. 
We also let $\bar{u}^i$ be the solution to
\begin{equation} \label{mean}
\begin{split}
\max_{u^i \in \mathbb{U}^i} \quad &\mathbb{E} [ \bar{J}_i^P [ u^i, 0, 0]  ]  \\
\mbox{subject to} \quad & 
dw_t = r_0 (\nu(t) - w_t) dt + \sigma_0 (t) dW_t^0 \\
&dx_t^i = f_i(x_t^i, u_t^i) dt.
\end{split}
\end{equation}
Suppose that $\mathbb{E}[ \bar{J}^P [ u, 0, 0]  ] > 0$, where $\bar{u} := (\bar{u}^1, \cdots, \bar{u}^n)$.
Set
\begin{equation} \nonumber
\rho := \frac{- \frac{1}{\theta} \log \mathbb{E}\left [\exp \left (-\theta \bar{J}^P [u^*, \gamma^*, \zeta^*]  \right )  \right]}{\mathbb{E} [\bar{J}^P [\bar{u}, 0, 0]]}.
\end{equation}
Then, the following suboptimality bound holds:
\begin{equation} \label{bound}
\begin{split}
&\rho \left (- \frac{1}{\theta} \log \mathbb{E}\left [\exp \left (-\theta \bar{J}^P [u^{\mbox{\tiny OPT}}, \gamma^{\mbox{\tiny OPT}}, \zeta^{\mbox{\tiny OPT}}]  \right )  \right] \right ) \leq - \frac{1}{\theta} \log \mathbb{E}\left [\exp \left (-\theta \bar{J}^P [u^*, \gamma^*, \zeta^*]  \right )  \right]
\end{split}
\end{equation}
for $\theta > 0$.
\end{proposition}

\begin{IEEEproof}
Using Jensen's inequality, we have
\begin{equation} \nonumber
\begin{split}
- \frac{1}{\theta} \log \mathbb{E}\left [\exp \left (-\theta \bar{J}^P [u^{\mbox{\tiny OPT}}, \gamma^{\mbox{\tiny OPT}}, \zeta^{\mbox{\tiny OPT}}]  \right )  \right]
&\leq - \frac{1}{\theta} \log \exp \left ( - \theta \mathbb{E} [\bar{J}^P [u^{\mbox{\tiny OPT}}, \gamma^{\mbox{\tiny OPT}}, \zeta^{\mbox{\tiny OPT}}]] \right )\\
&=
 \mathbb{E} [\bar{J}^P [u^{\mbox{\tiny OPT}}, \gamma^{\mbox{\tiny OPT}}, \zeta^{\mbox{\tiny OPT}}]]. 
\end{split}
\end{equation}
When the principal is risk-neutral, 
the principal's payoff is independent of
$(\gamma, \zeta)$ and setting $(\gamma, \zeta) = (0,0)$ always satisfies the risk-limiting condition.
From this observation, we claim that $(\bar{u}, 0, 0)$ solves the problem \eqref{opt} when $\theta = 0$, i.e., the objective function is replaced with $\mathbb{E} [ \bar{J}^P [u, \gamma, \zeta]]$.
We first note that $(\bar{u}, 0, 0)$ satisfies all the constraints in \eqref{opt}. 
Suppose that $(\bar{u}, 0, 0)$ does not solve \eqref{opt} and choose a solution, $(\hat{u}, \hat{\gamma}, \hat{\zeta})$, of \eqref{opt}.
Because $\mathbb{E} [ \bar{J}_i^P [ \hat{u}^i, \hat{\gamma}, \hat{\zeta}] ] = \mathbb{E} [ \bar{J}_i^P [ \hat{u}^i, 0, 0] ]$ and $\hat{u}$ satisfies the constraint of \eqref{mean}, the following inequality holds:
\begin{equation} \nonumber
\mathbb{E} [ \bar{J}_i^P [ \hat{u}^i, \hat{\gamma}^i, \hat{\zeta}^i] ] = \mathbb{E} [ \bar{J}_i^P [ \hat{u}^i, 0, 0] ] \leq \mathbb{E} [ \bar{J}_i^P [ \bar{u}^i, 0, 0] ].
\end{equation}
Therefore, we have
\begin{equation} \nonumber
\begin{split}
\mathbb{E} [ \bar{J}^P [ \hat{u}, \hat{\gamma}, \hat{\zeta}] ] &= \sum_{i=1}^n \mathbb{E} [ \bar{J}_i^P [ \hat{u}^i, \hat{\gamma}^i, \hat{\zeta}^i] ]\\
& \leq \sum_{i=1}^n \mathbb{E} [ \bar{J}_i^P [ \bar{u}^i, 0, 0] ]
=\mathbb{E} [ \bar{J}^P [ \bar{u}, 0, 0] ].
\end{split}
\end{equation}
This inequality is contradictory to the fact that $(\bar{u}, 0, 0)$ does not solve \eqref{opt}. 
Therefore, 
\begin{equation} \nonumber
\mathbb{E} [\bar{J}^P [u^{\mbox{\tiny OPT}}, \gamma^{\mbox{\tiny OPT}}, \zeta^{\mbox{\tiny OPT}}]] 
\leq \mathbb{E} [\bar{J}^P [\bar{u}, 0,0] ].
\end{equation}
As a result, the suboptimality bound \eqref{bound} holds.
\end{IEEEproof}
Note that the suboptimality bound can be computed by solving \eqref{dec} and \eqref{mean} for each $i$ while it is not feasible to directly solve \eqref{opt} for $n > 1$.
This proposition implies that the proposed decomposition tends to be exact as the coefficient $\theta$ of the principal's risk aversion goes to zero because $\rho \to 1$ as $\theta \to 0$.
Furthermore, the approximate contract $(C^*, u^*)$ satisfies the participation-payoff condition \eqref{individual} and the risk-limiting condition \eqref{risk_limiting}.

%
%
%
%
Due to this decomposition, the contract design for agent $i$ only requires the state space, $\mathbb{R}^3$, of $(w_t, x_t^i, y_t^i)$ rather than the full joint state space, $\mathbb{R}^{3n+1}$, of $(w_t, x_t, v_t, y_t)$.
Therefore, the computational complexity of designing a risk-limiting contract for an agent is independent of the total number of agents.
The decomposed problem for agent $i$ is solved via dynamic programming over the reduced state space, $\mathbb{R}^3$, of $(w_t, x_t^i, y_t^i)$ as follows.
We set the feasible set of control as
\begin{equation}\nonumber
\Omega^i := \{(u^i, \gamma^i, \zeta^i) \in \mathbb{U}^i \times \Gamma^i \times \Gamma^i \: | \: y_T^i \geq 0 \mbox{ a.s.} \}.
\end{equation}
To synthesize a risk-limiting dynamic contract for agent $i$, we first define the value function of \eqref{dec} associated with agent $i$ as
\begin{equation} \label{value}
\begin{split}
&\phi_i(\bm{w}, \bm{x}_i,\bm{y}_i,t) := \max_{(u^i, \gamma^i, \zeta^i) \in \Omega^i} \;  -\frac{1}{\theta} \log \mathbb{E}_{\bm{w}, \bm{x}_i, \bm{y}_i,t}
\\
&\qquad \quad\left [ \exp \left (- \theta \left ( \int_t^T R_i(w_s, x_s^i, u_s^i) ds   + \int_t^T G_i(w_s, \gamma_s^i) dW_s^{(i)}  - b_i  \right ) \right )  \right ],
\end{split}
\end{equation}
where 
$\mathbb{E}_{\bm{w}, \bm{x}_i, \bm{y}_i, t} [A]$ denotes the expectation of $A$ conditioned on $(w_t, x_t^i, y_t^i) = (\bm{w}, \bm{x}_i, \bm{y}_i)$, and 
\begin{equation} \nonumber
\begin{split}
R_i(\bm{w}, \bm{x}_i, \bm{u}) &:= r_i^P (\bm{w}, \bm{x}_i, \bm{u}) + r_i^A (\bm{x}_i, \bm{u}),\\
G_i(t, \bm{w}, \bm{\gamma}) &:= \begin{bmatrix}
-\bm{\gamma}_1 & -\bm{\gamma}_2 +  \sigma_i^P (\bm{w}) + \sigma_i^A(t))
\end{bmatrix}. 
\end{split}
\end{equation}
To handle the constraint $y_T^i \geq 0$ a.s., which is often called the \emph{stochastic target constraint}, we use the Hamilton-Jacobi-Bellman (HJB) characterization proposed in \cite{Bouchard2010}.
This characterization converts the target constraint into a `classical' state constraint using the \emph{geometric dynamic programming principle} \cite{Soner2002}.
The reformulated constraint  is embedded in an auxiliary value function. 
This auxiliary value function is a viscosity solution of an HJB equation.
Applying the dynamic programming principle on the original value function \eqref{value}, one can derive a constrained-HJB equation in which the stochastic target constraint is reformulated as the constraints on the auxiliary value function and the control. 
In our case, the auxiliary value function is a zero function. 
Let 
\begin{equation}\nonumber
U_i(\bm{y}_i) := \{(\bm{u}, \bm{\gamma}, \bm{\zeta}) \in \mathcal{U}^i \times \mathbb{R}^2 \times \mathbb{R}^2 \: | \: \bm{\gamma} = \bm{\zeta} = 0 \mbox{ if }\bm{y}_i \leq 0 \}.
\end{equation}
Then, the stochastic target constraint is simply incorporated into the following constrained-HJB equation:
\begin{equation}\nonumber
\begin{split}
&\frac{\partial \phi_i}{\partial t}  +\max_{\substack{(\bm{u}, \bm{\gamma}, \bm{\zeta})\\ \in U_i(\bm{y}_i)}}
\left \{
(F_i(\bm{w}, \bm{x}_i, \bm{u}, \bm{\gamma}) - \theta \Sigma(\bm{\zeta} )G_i(\bm{w}, \bm{\gamma})^\top  )^\top D \phi_i \right.\\
&\qquad \quad \left. +R_i(\bm{w}, \bm{x}_i, \bm{u}) 
-\frac{\theta}{2} \|G_i(\bm{w}, \bm{\gamma})  \|^2
-\frac{\theta}{2} \| \Sigma(\bm{\zeta} )^\top D\phi_i \|^2
 + \frac{1}{2}\mbox{tr} (\Sigma(\bm{\zeta} ) \Sigma(\bm{\zeta} )^\top D^2 \phi_i) 
\right\} = 0,
\\
&\phi_i (\bm{w}, \bm{x}_i, \bm{y}_i, T) = -b_i,
\end{split}
\end{equation}
whose viscosity solution corresponds to the value function \eqref{value} \cite{Crandall1992, Fleming2006, Bouchard2010},
where
\begin{equation} \nonumber
\begin{split}
F_i(t,\bm{w}, \bm{x}_i, \bm{u}, \bm{\gamma}) := \begin{bmatrix}
r_0 (\nu (t) - \bm{w})\\
f_i (\bm{x}_i, \bm{u}) \\
-\| \bm{\gamma}\|^2
\end{bmatrix}, \quad
\Sigma(t, \bm{\zeta}) := \begin{bmatrix}
\sigma_0 (t) & 0 \\
0 & 0\\
\bm{\zeta}_1 & \bm{\zeta}_2
\end{bmatrix}.
\end{split}
\end{equation}
In general, an analytic solution of the HJB equation is difficult to find. 
Therefore, we grid up the state space, which corresponds to the domain of the PDE, and we numerically evaluate the solution at the grid points using convergent schemes, e.g., \cite{Barles1991, Kushner2001}. 
After solving the HJB equation, we can use the value function to compute an optimal compensation scheme, $C^*$, and an optimal control strategy, $u^*$.
Set $(w_0, x_0^{*i}, y_0^{*i}) = (\ln {\lambda^0}, x^{0i}, S_i)$ for $i=1, \cdots, n$. Given the process $\bold{x}_s^{*i} := (w_s^{*}, x_s^{*i}, y_s^{*i}) \in \mathbb{R}^3$ for $s \in [0,t]$, we can determine an optimal control strategy as
\begin{equation}\nonumber
\begin{split}
(u_t^{*i}, \gamma_t^{*i}, \zeta_t^{*i}) = \arg \max_{(\bm{u}, \bm{\gamma}, \bm{\zeta})\in U_i(\bm{y}_i)}
&\left \{
(F_i(w_t, x_t^{*i}, \bm{u}, \bm{\gamma}) - \theta \Sigma(\bm{\zeta} ) G_i(w_t, \bm{\gamma})^\top )^\top D \phi_i (\bold{x}_t^{*i}, t) \right.\\
&+R_i(w_t, x_t^{*i}, \bm{u}) 
-\frac{\theta}{2}\|G_i(w_t, \bm{\gamma})  \|^2 -\frac{\theta}{2} \| \Sigma(\bm{\zeta} )^\top D\phi_i (\bold{x}_t^{*i}, t) \|^2
\\
&\left. 
+ \frac{1}{2}\mbox{tr} \left (\Sigma(\bm{\zeta} ) \Sigma(\bm{\zeta} )^\top D^2 \phi_i (\bold{x}_t^{*i}, t) \right ) 
\right\}
\end{split}
\end{equation}
 for $i = 1, \cdots, n$ and
for $t \in [0,T]$.
Note that $v_t^{*i}$ and $y_t^{*i}$  can be computed by integrating the SDEs \eqref{v_s} and \eqref{sde2}
over time with the control $(u_s^{*i}, \gamma_s^{*i}, \zeta_s^{*i})$ for $s \in [0,t)$, respectively.
Then, an optimal compensation scheme can be obtained as
$C^* = v_T^*$
which is proposed in Theorem \ref{optimal}.
A more detailed discussion regarding how to synthesize an optimal control using a viscosity solution of an associated HJB equation can be found in \cite{Bardi1997} even when the viscosity solution is not differentiable.

\begin{figure}[tb] 
\begin{center}
\includegraphics[width =2.9in]{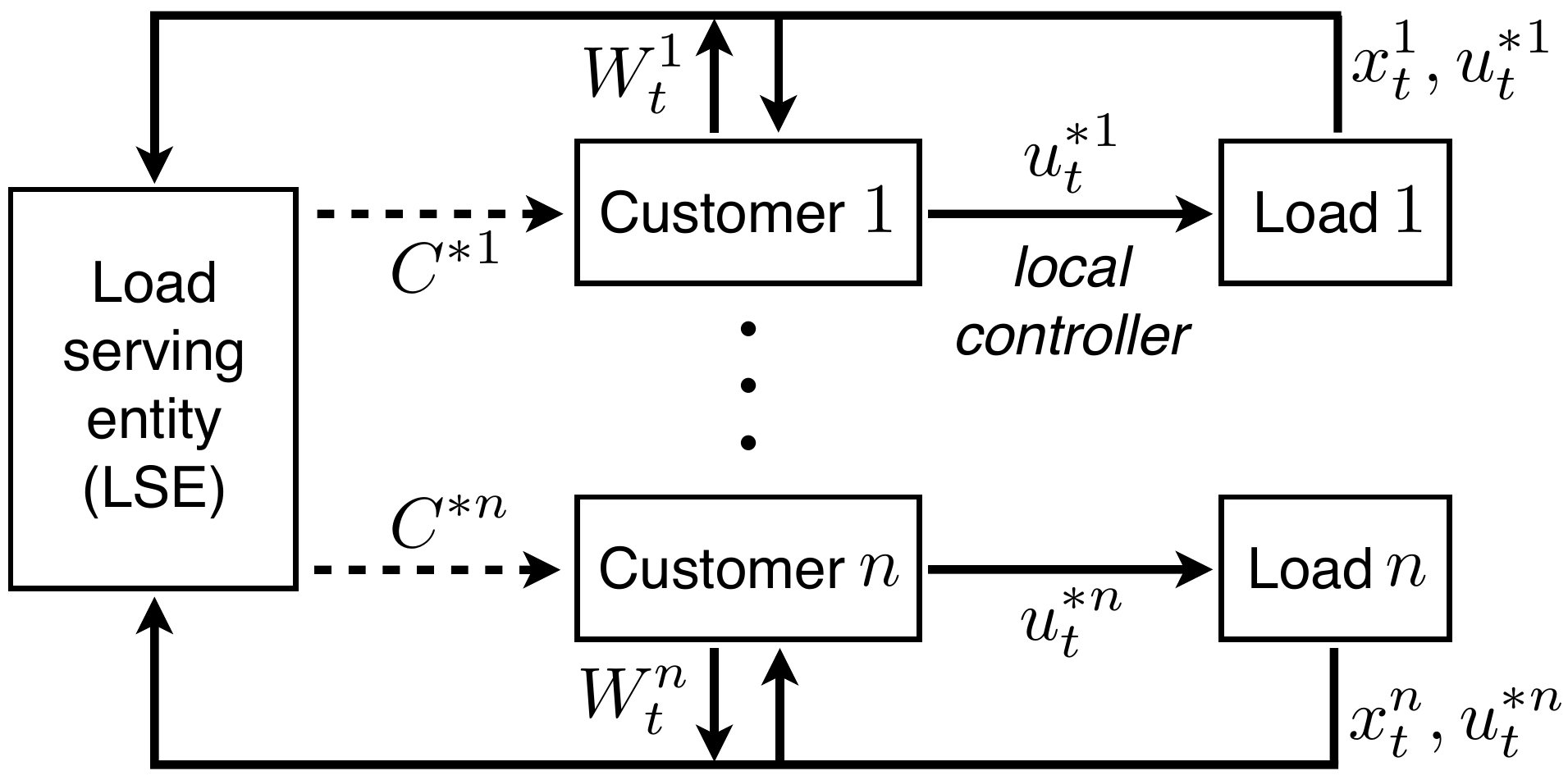}
\caption{Implementation of the proposed contracts: the controls of loads can be decentralized with a broadcast of price (LMP) information, while a centralized monitoring is required. The compensations are provided at the end of the contract period.
}
 \label{fig:dlc}
 \end{center}
\end{figure}

The synthesized optimal control is written into the contract and the customer must follow the control strategy if he or she enters into the contract. 
Note that the control for one load is independent of that for another load. 
Furthermore, the optimal control for a load is given as state-feedback, where the state variables only require the energy price in the real-time market and the local information of the load.
Therefore, the proposed control can be decentralized with a broadcast of the price information, i.e., the local controller in which the optimal control strategy is programmed is sufficient for the implementation of the contract as depicted in Fig. \ref{fig:dlc}.
On the other hand, the load-serving entity still needs to monitor the state, the control and the forecast error for each customer to ensure that each customer follows the optimal control strategy written in the contract. 
The total power consumption of each customer monitored by a smart meter can be used to compute the forecast error. In addition, the sensors for loads such as thermostats provide the state and control information. 
The monitored information could be transferred to the load-serving entity through a one-way data connection such as the Internet.
The information gathered by the monitoring is also used to
 compute the optimal compensation provided to each customer at the end of the contract period.

\section{Application to Direct Load Control for Financial Risk Management} \label{application}

In this section, we apply the proposed risk-limiting dynamic contracts to direct load control. 
The performance and usefulness of the novel direct load control program for financial risk management are demonstrated using the data of LMPs in the ERCOT and the electric energy consumption of customers in Austin, Texas.

\subsection{Data and Setting}

\begin{figure}[tb] 
\begin{center}
\includegraphics[width =3.5in]{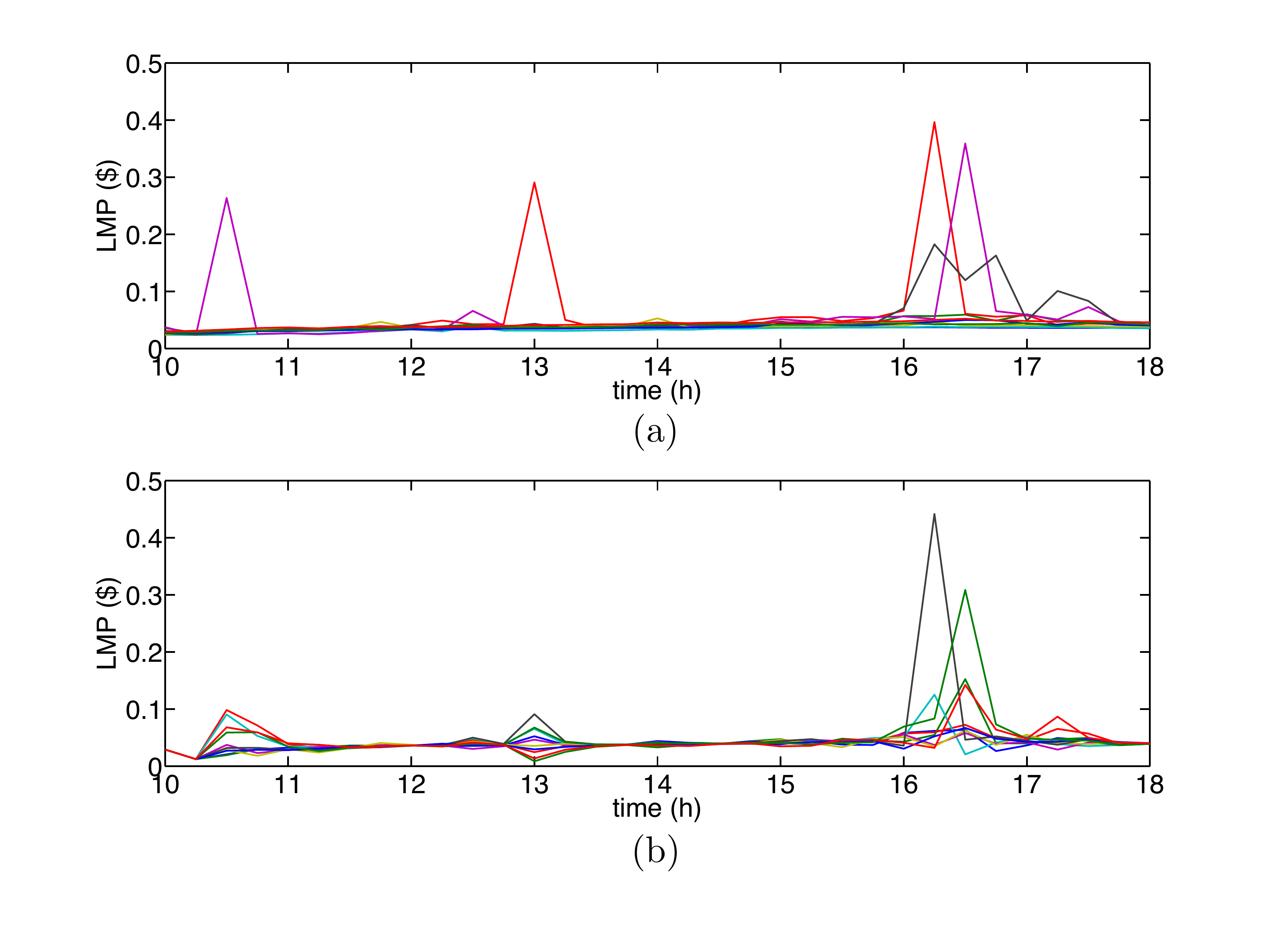}
\caption{(a) Locational marginal price (LMP) in Austin, Texas, from July 1, 2013 to July 10, 2013. (b) Ten sampled trajectories of LMP, $\{ \lambda_t\}_{0\leq t \leq T}$, generated by our identified price model \eqref{mr}. 
}
 \label{fig:lmp}
 \end{center}
\end{figure}  

We consider a scenario in which each customer provides one air conditioner for the proposed direct load control program.
We use the ETP model \eqref{indoor_ex} for the customer's indoor temperature dynamics given in Example \ref{ex_tcl}.
The set of feasible control values is chosen as $\mathcal{U}^i := \{ 0, 2\}$, assuming that customer $i$'s air conditioner consumes $0$kW in its OFF state and $2$kW in its ON state.
The model parameters are chosen as $\alpha_i = 0.1$ and $\kappa_i = 1.5$, which are calculated based on the Residential module user's guide from GridLAB-D and are physically reasonable \cite{gridlab}.
Customer $i$'s comfort level is chosen as \eqref{comfort_tcl} in Example \ref{ex_comfort}, with $\omega_i = 0.15$.
We set the customer's desirable indoor temperature range, $[\underline{\Theta}, \overline{\Theta}]$, as $[20^\circ\mbox{C}, 22^\circ\mbox{C}]$.
We choose the contract period as $[10\mbox{h}, 18\mbox{h}]$.

We use customers' electric energy consumption data  in Austin, Texas \cite{wikienergy} to estimate the load profile $l_i (t)$ and the diffusion coefficient $\tilde{\sigma}_i (t)$ in \eqref{demand}. The load $l_i(t)$ is chosen as the mean value of the customer $i$'s power consumption at $t$ other than the air conditioner.
For the estimation of $\tilde{\sigma}_i$, we apply the Kalman filter \cite{Kalman1961, Kristensen2004} over the data set for the summer period, from June to September 2013, assuming that customer $i$'s energy consumption profile other than the air conditioner for one day in the period represents one sampled trajectory.  
We then scale the estimated diffusion coefficient by a constant factor such that $\int_0^T \tilde{\sigma}_i(t)^2 dt$ is equal to the variance of the energy consumption data.
This scaling guarantees that $\mbox{Var} [J_i^A[0,0]]$ is equal to the variance of the customer's energy cost.

The price model \eqref{mr} is identified  using the ERCOT LMP data  at the settlement point, AUSTIN PLANT, from July 1, 2013 to July 10, 2013 \cite{ercot}. We estimate the parameters, $r_0$, $\nu(t)$ and $\sigma_0(t)$, by applying the Kalman filter on the transformed linear model \eqref{mr2}. 
The LMP data and samples of the price profile generated by the identified model are shown in Fig. \ref{fig:lmp}.

\begin{figure}[tb] 
\begin{center}
\includegraphics[width =3.4in]{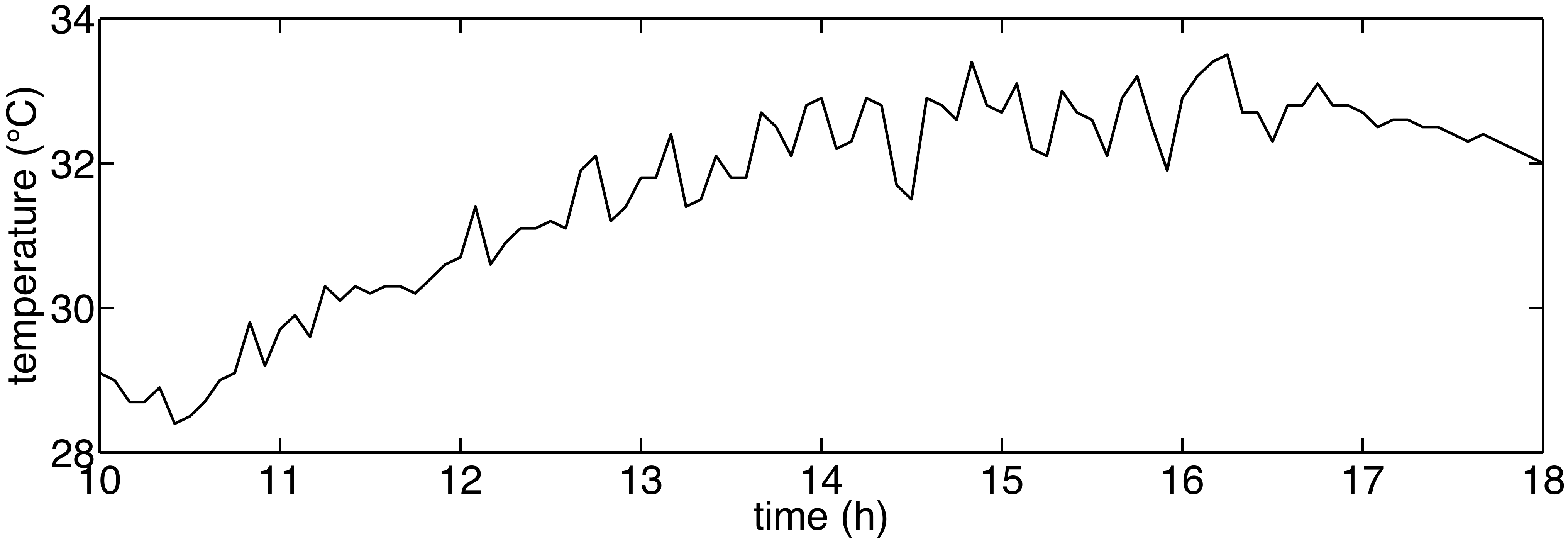}
\caption{Outdoor air temperature, $\Theta(t)$, in Austin, Texas, on July 5, 2013.
}
 \label{fig:temp_austin}
 \end{center}
\end{figure}

We use the NOAA Quality Controlled Local Climatological Data in Austin, Texas for the outdoor temperature profile \cite{noaa}. The outdoor temperature, $\Theta(t)$, at time $t \in [10\mbox{h}, 18\mbox{h}]$ is chosen as the temperature on July 5, 2013, at time $t \in [10\mbox{h}, 18\mbox{h}]$ and is shown in Fig. \ref{fig:temp_austin}.

\subsection{Comparison to Optimal Load Control by Customers}

Suppose that customer $i$ does not participate in the direct load control program and has the following payoff:
\begin{equation}
\hat{J}_i^A [u^i] := \int_0^T r_i^A(x_t^i, u_t^i) dt + \int_0^T \sigma_i^A (t) dW_t^i,
\end{equation}
where $r_i^A$ and $\sigma_i^A$ are given by \eqref{a_ex}.
The solution of the following optimal control problem maximizes customer $i$'s expected payoff:
\begin{equation} \label{det_opt}
\begin{split}
\max_{u^i \in \mathbb{U}^i} \quad &\mathbb{E} [ \hat{J}_i^A[u^i] ]:=  \int_0^T r_i^A(x_t^i, u_t^i) dt \\
\mbox{subject to} \quad & dx_t^i = f_i (x_t^i,u_t^i) dt.
\end{split}
\end{equation}
The optimal control can be obtained using the viscosity solution of the following HJB equation \cite{Bardi1997}:
\begin{equation} \nonumber
\begin{split}
&\frac{\partial \hat{\phi}_i(\bm{x}_i, t)}{\partial t} + \max_{\bm{u} \in \mathcal{U}^i} \{
f(\bm{x}_i, \bm{u}) D_{\bm{x}_i} \hat{\phi}_i (\bm{x}_i, t) 
+ r_i^A(\bm{x}_i, \bm{u}) 
\} = 0,\\
&\hat{\phi}_i (\bm{x}_i, T) = 0.
\end{split}
\end{equation}
Let $\hat{u}^{*i}$ be an optimal control, i.e., a solution to \eqref{det_opt}. We let $\bar{b}_i := \mathbb{E} [\hat{J}_i^A[\hat{u}^{*i}]]$ and $\bar{S}_i := \mbox{Var} [\hat{J}_i^A[\hat{u}^{*i}]] = \int_0^T \sigma_i^A(t)^2 dt$ be the nominal expected payoff and risk of customer $i$, respectively.

\begin{figure}[tb] 
\begin{center}
\includegraphics[width =3.5in]{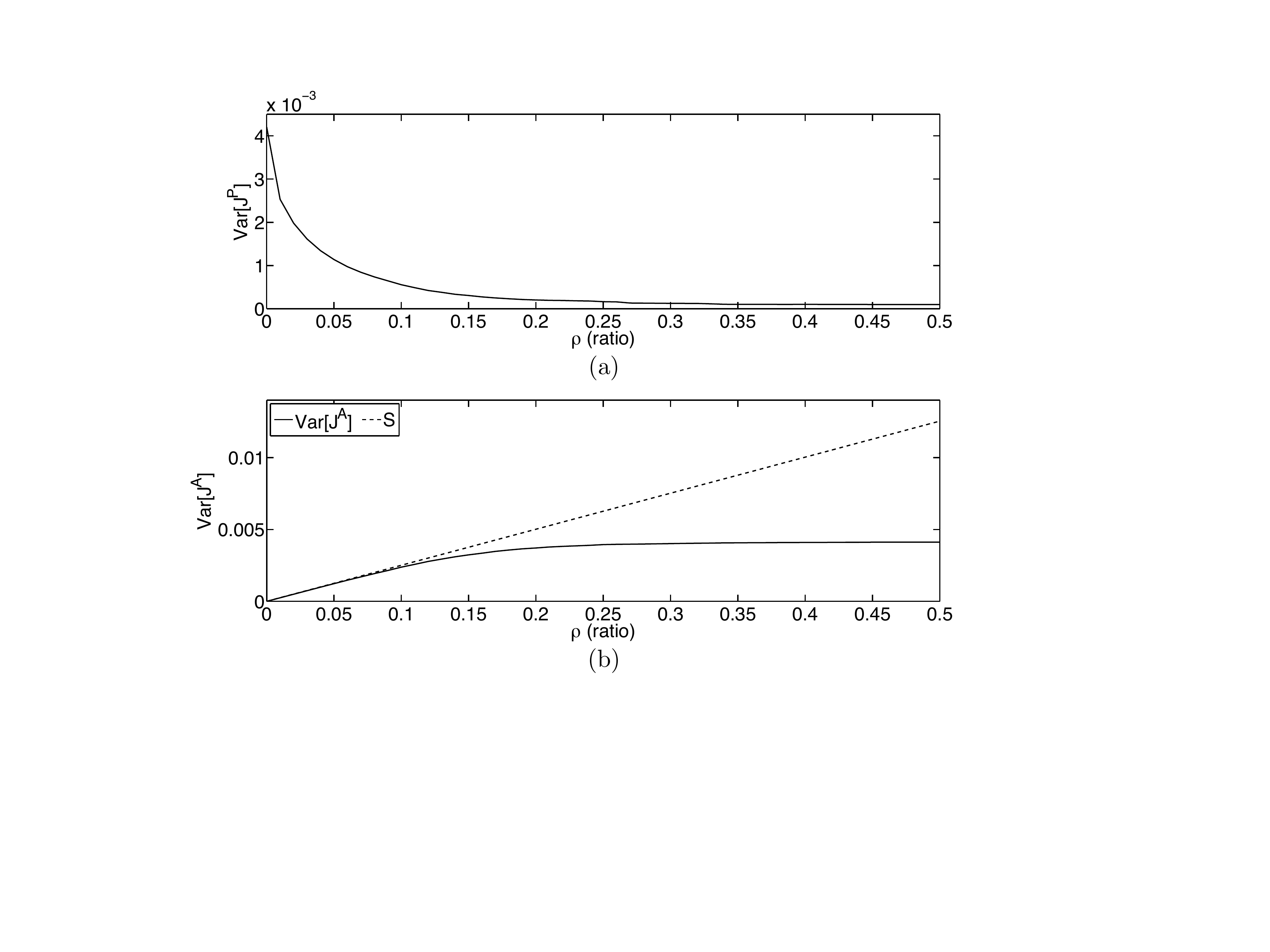}
\caption{The simulation results with the contract: (a) variance of the load-serving entity's payoff.
(b) variance of the customer's payoff and
the risk share $S_i = \rho\bar{S}_i$.
}
 \label{fig:result1e3}
 \end{center}
\end{figure}

We now compare the performance of the proposed contracts to that of this optimal load control without a contract.
We choose the customer's electricity price as the flat price, $\mu_i \equiv \bar{\mu} = \$ 0.11$, specified in Austin Energy's electricity tariff for summer \cite{austin}.
In the absence of a contract, the mean and variance of customer $i$'s payoff are $\bar{b}_i$ and $\bar{S}_i$, respectively. 
We set $b_i = \bar{b}_i$ and $S_i = \rho \bar{S}_i$ and vary $\rho$ from $0$ to $0.3$. 
If customer $i$ enters into the contract, which is the solution of \eqref{dec}, then the mean value of the customer's payoff is guaranteed to be greater than or equal to $\bar{b}_i$, and the variance of the customer's payoff is guaranteed to be less than or equal to $\rho\bar{S}_i$. 

The  results of numerical experiments presented in Fig. \ref{fig:result1e3} verify the performance of the proposed contract. 
The coefficient of load-serving entity's risk aversion is chosen as $\theta = 10^{-2}$.
As shown in Fig. \ref{fig:result1e3} (a), The variance of the load-serving entity's payoff decreases as the ratio $\rho$ of the amount of the risk that the customer is willing to bear to the customer's nominal risk increases.
On the other hand,
the variance of the customer's payoff (blue) increases as the ratio $\rho$ of the amount of the risk that the customer is willing to bear to the customer's nominal risk increases, as shown in Fig. \ref{fig:result1e3} (b).
More importantly, it is less than or equal to the risk share $S_i = \rho\bar{S}_i$ (red).
Therefore, we confirm that the risk-limiting condition is satisfied.

When the customer does not enter into in the contract, the variance of the load-serving entity's payoff is 
\begin{equation}\nonumber
\begin{split}
&\mbox{Var}[\hat{J}_i^P[\hat{u}^{*i}]] \\
&:= \mbox{Var} \left [\int_0^T r_i^P(w_t, \hat{x}_t^{*i}, \hat{u}_t^{*i}) dt + \int_0^T \sigma_i^P(w_t) dW_t^i \right] \\
& = 0.0108.
\end{split}
\end{equation}
By comparing this variance with the variance of the load-serving entity's payoff when the contract is executed (Fig. \ref{fig:result1e3} (a)),
we note that the contract reduces the load-serving entity's risk by more than 50\% even when the customer is extremely risk-averse, i.e., $S_i = 0$.
If the customer chooses $S_i \geq 0.2 \bar{S}_i$ in the contract, the load-serving entity's risk is decreased by more than 95\%.

The mean values of the load-serving entity's payoff with and without the contract are given by
\begin{equation} \nonumber
\mathbb{E} [ J_i^P [C^{*i}, u^{*i} ]] = 1.324, \quad \mathbb{E} [\hat{J}_i^P[\hat{u}^{*i}]] = 1.297,
\end{equation}
respectively. 
Therefore, the load-serving entity can pay $\$ 0.027$ more for the customer without reducing its mean payoff if the customer enters into a contract.
In other words, the load-serving entity can incentivize the customer to enter into the contract by increasing the customer's expected payoff by $\$ 0.027$.


\begin{figure}[tb] 
\begin{center}
\includegraphics[width =3.5in]{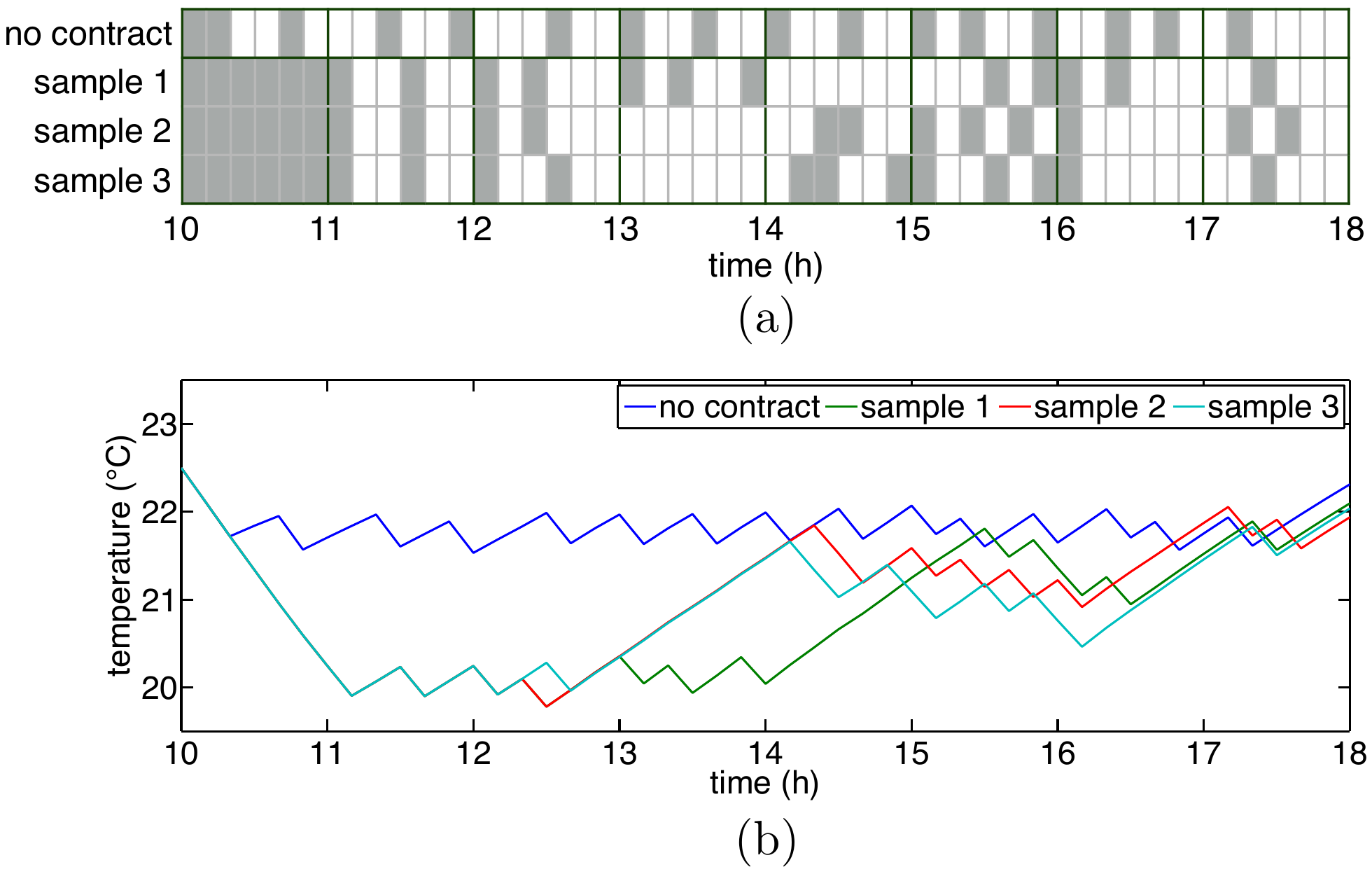}
\caption{The simulation results without and with the contract (three samples are presented in the case with the contract): (a) control (gray: ON, white: OFF); and
(b) indoor temperature.
}
 \label{fig:sample_result1e3}
 \end{center}
\end{figure}  

Fig. \ref{fig:sample_result1e3} shows the effect of the contract on the control and the indoor temperature. 
In this set of experiments, we set $S_i = 0.1 \bar{S}_i$ and $\theta = 10^{-2}$.
When the customer does not enter into the contract, the customer regularly turns on and off the air conditioner such that the indoor temperature is kept near 22$^\circ$C. 
This is because the customer wants to save the energy cost by only taking into account the energy price $\mu_i$, which is fixed for all time.
On the other hand, if the customer enters into the contract, the room is pre-cooled before 12:30pm when the LMP $\lambda_t$ is low and not volatile. 
Another pre-cooling interval is from 3pm to 4pm. This pre-cooling allows to save energy purchase in the real-time market from 4pm to 5pm when the LMP is highly volatile. 
In other words, the contract properly manages the price risk in the spot market.

%

\subsection{Validation of the Brownian Motion Model Using Data}

The energy consumption process of a customer is modeled by the SDE \eqref{demand}.
In practice, the load forecast error may not be exactly captured by the diffusion term, $\tilde{\sigma}_i (t) dW_t^i$, with a standard Brownian motion.
We test the robustness of the proposed contract method with respect to the deviation of the demand forecast errors in the data from the Brownian motion model.
More specifically, we execute the optimal contract synthesized using the Brownian motion model over the data.
Then, we compute how much the load-serving entity's and the customer's resulting payoffs differ from their optimal payoffs obtained under the Brownian motion assumption. 
Setting $S_i = \rho \bar{S}_i$,
we averaged the percentage deviations over $\rho = [0, 1]$.
The average deviations in the mean of the load-serving entity's and the customer's payoffs are  $0.010\%$ and $0.012\%$, respectively. Furthermore, the risk-limiting condition is not violated for $\rho > 0.14$; furthermore, for $\rho \leq 0.14$, the condition is never violated by more than 12\%.
This preliminary test of the proposed contract framework with respect to errors in the Brownian motion model suggests that the framework is robust to load model errors, however
more data is need for accurate approximations of the mean and variance values.
Further experiments will be performed in the future to rigorously test the validity of Brownian motion model in the proposed contracts.

\section{Conclusions}

We proposed a new continuous-time dynamic contract framework that has a risk-limiting capability.
The key feature of the proposed contract is that the variance of the agent's payoff is bounded by a threshold specified in the contract.
This feature enables the contract framework, by combining with direct load control, to provide financial risk management solutions for real-time electricity markets.
To obtain a globally optimal contract, a dynamic programming-based method is developed.
Difficulty arises in dealing with the constraints on the mean and the variance of the agent's payoff. 
We resolve this issue by introducing two new dynamical systems that have intuitive meanings. 
We also proposed an approximate decomposition of
the contract design problem for $n$ agents into $n$ low-dimensional problems for each agent. 
The approximate contract obtained using the proposed decomposition has a guaranteed suboptimality bound.
This decomposability allows the direct load control program based on this contract framework to handle a large number of customers without any scalability issue.
Using the data of the ERCOT LMPs and electric energy consumption by customers in Austin, Texas, we perform numerical experiments to validate the performance and usefulness of the proposed contracts.


\section*{Acknowledgement}

The authors would like to thank Professor Lawrence Craig Evans for helpful discussions on dynamic contract theory and PDE-based approaches for it.

\bibliographystyle{IEEEtran}

\bibliography{RLDC_DLC_onecolumn}
\vfill\eject

\end{document}